\documentclass[11pt]{amsart}
\usepackage{rotating}
\usepackage{mathcomp,amscd}
\usepackage{amssymb}
\usepackage{hyperref}
\usepackage{times}
\usepackage{enumerate}
\usepackage{color}

\newcommand{\GL}{\mathrm{GL}}

\newcommand{\wt}{\mathrm{wt}}
\newcommand{\noJ}{\# \mathcal{J}}

\newcommand{\diag}{\mathrm{diag}}

\newcommand{\subspace}{Q}

\usepackage[OT2,T1]{fontenc}
\DeclareSymbolFont{cyrletters}{OT2}{wncyr}{m}{n}
\DeclareFontFamily{OT1}{rsfs}{}
\DeclareFontEncoding{OT2}{}{} 
  
     \DeclareFontShape{OT1}{rsfs}{n}{it}{<-> rsfs10}{}
\DeclareMathAlphabet{\mathscr}{OT1}{rsfs}{n}{it}

\newcommand{\C}{{\Bbb C}}
\newcommand{\Z}{{\Bbb Z}}

\newcommand{\R}{{\Bbb R}}

\newcommand{\vol}{\mathrm{vol}}

\newcommand{\SL}{\mathrm{SL}}

\numberwithin{equation}{section}
\numberwithin{table}{section}
\numberwithin{figure}{section}
\newtheorem{Theorem}{Theorem}[section]

\newtheorem{Lemma}{Lemma}[section]

\allowdisplaybreaks
\begin{document}
   \title{The behavior of random reduced bases} 
    \author{Seungki Kim and Akshay Venkatesh}

\begin{abstract}     
We prove that the number of  Siegel-reduced bases  for a randomly chosen $n$-dimensional lattice 
becomes, for $n \rightarrow \infty$,  tightly concentrated around its mean. We also show that most  reduced bases behave as in the worst-case analysis of lattice reduction. 
    Comparing with experiment, these results  suggest   that most reduced bases     will, in fact,  ``very rarely'' occur as an output of lattice reduction. 
The concentration result is based on an analysis of the spectral theory of Eisenstein series and uses (probably in a removable way) the Riemann hypothesis. 

   \end{abstract}
 \maketitle
 \section{Introduction}
 For us, a lattice $L \subset \mathbf{R}^n$ is the set $\Z. \mathcal{B}$ of all linear combinations  of a basis $\mathcal{B} = \{\mathbf{x}_1, \dots, \mathbf{x}_n\}$
 for $\R^n$; we say, then, that $\mathcal{B}$ is a basis for $L$.  The volume $\vol(L)$ of the lattice $L$ is the determinant of the matrix with rows $\mathbf{x}_i$.
 In what follows, we assume $\vol(L)=1$. 
 
Define 
 $  \mathbf{x}_i^*$ to be   the projection of $\mathbf{x}_i$ to the orthogonal complement of $\langle \mathbf{x}_{i+1},\dots, \mathbf{x}_n \rangle$. 
 We say that the basis $\mathcal{B} = (\mathbf{x}_1, \dots, \mathbf{x}_n)$ is Siegel-reduced with parameter $T$ if
 the following conditions hold: 
 \begin{itemize}
 \item
 $\|\mathbf{x}_i^*\| \geq T^{-1} \| \mathbf{x}_{i+1}^*\|$, and
 \item If we write $\mathbf{x}_i = \mathbf{x}_i^* + \sum_{j > i} n_{ji} \mathbf{x}_j^*$, then all $|n_{ji}| \leq \frac{1}{2}$. 
 \end{itemize}
  Since we suppose $\vol(L)=1$, we obtain\footnote{Our indexing of the basis corresponds to the standard numbering of roots for $\SL_n$, and is unfortunately opposite to that usually used in analysis of LLL.}
 \begin{equation} \label{worst-case} \|\mathbf{x}_n\| \leq T^{\frac{n-1}{2}} .\end{equation}
   
 The LLL algorithm produces (in polynomial time)  \footnote{ In fact, they satisfy a slightly stronger reduction condition. We ignore the difference for the purpose of this introduction. Our expectation is that very similar theorems hold in both cases, but we don't know how to prove our main result for the LLL reduction condition.}
 Siegel-reduced bases with parameter $T$ for any lattice $L$ and any $T>T_0 := 2/\sqrt{3}$.  In particular,
 it produces a  ``relatively short'' vector $\mathbf{x}_n$, which is guaranteed to satisfy \eqref{worst-case}.
  
 In practice, the situation is even better:  Nguyen and Stehl{\'e} \cite{NS} have investigated in detail the experimental behavior of the LLL algorithm
 and observed that 
 it ``typically'' produces a basis with $\|\mathbf{x}_n\| \approx (1.02)^n$. By comparison, $\sqrt{T_0} \approx 1.075 \dots$; said differently,
the typical quality of an output basis of LLL is very much better than the worst-case bound  \eqref{worst-case} for reduced bases. 
 
The main point of this paper is to observe that the output of LLL is not just better than the worst-case bound for reduced bases, but also better than
the average bound for reduced bases.  Recall \cite{Siegel} that there is a unique probability measure $\mu_n$ on the space of covolume $1$ lattices
which is invariant by linear transformations; thus there is a notion of random lattice.   The following gives a flavor of what's proven: 

\begin{quote} 
(*) If  we   first
choose a $\mu_n$-random lattice  $L$, and then choose a basis $\mathcal{B}$ uniformly and randomly from the finite set  $\{\mathcal{B}_1, \dots, \mathcal{B}_r\}$ of  Siegel-reduced bases for $L$, 
 we have $\frac{\|\mathbf{x}_n\|}{T^{(n-1)/2}}  > 0.999$ with probability approaching $1$ as $n \rightarrow \infty$. 
\end{quote}

This result says that typical reduced bases behave just as badly as \eqref{worst-case}. It is derived from the more precise theorem below. 
Thus,  {\em the good properties
 of LLL are not merely a function of the properties of random reduced bases; the LLL algorithm itself ``selects'' good bases.}
This suggests, for example,  that there should be a very large number of ``dark'' reduced bases which are practically never selected by the LLL algorithm. 
It also suggests the importance of the following (not quite well-defined)

{\bf Problem.}  Determine a reliable heuristic
that, given a reduced basis $\mathcal{B}$, predicts how frequently it occurs as the output of LLL-reduction 
if we choose ``random'' input bases for the lattice $\Z.\mathcal{B}$.

 These questions have  been studied numerically  in the PhD thesis \cite{KimThesis} of the second-author; that work gives some evidence for the ``dark'' reduced basis phenomenon, and suggests that 
 that the likelihood of a reduced basis $(\mathbf{x}_i)_{1 \leq i \leq n}$ to be chosen by LLL   is inversely related to the ``energy'' $ \prod_{i=1}^{n} \| \mathbf{x}_i \wedge  \mathbf{x}_{i+1} \dots \wedge \mathbf{x}_n \| $.

To state the main theorem, we set up some notation. Let $T > 1$. Let
$$\mathfrak{S}_n   = \{ \mathcal{B}   \in (\mathbf{R}^n)^n: \mathcal{B} \mbox{ is Siegel-reduced with parameter $T$, }  \det(\Z.\mathcal{B}) = 1 \} $$be 
the set of  {Siegel}-reduced bases with parameter $T$ for lattices of volume $1$; here $\Z. \mathcal{B}$ denotes the lattice spanned by $\mathcal{B}$. Let $\mathcal{L}_n$ be the set
of lattices of determinant $1$. 
The natural map $$\pi: \mathfrak{S}_n \rightarrow \mathcal{L}_n, \ \ \pi(\mathcal{B}) = \Z.\mathcal{B}$$
has finite fibers. 
We equip $\mathcal{L}_n$ and $\mathfrak{S}_n$ with the probability measures invariant by $\SL_n(\R)$; then $\int_{\mathfrak{S}_n} f  = \int_{L \in \mathcal{L}_n} \sum_{\pi^{-1} L} f $. 
Let $N(L)$ be the  size of the fiber above $L \in \mathcal{L}_n$, i.e., the number of Siegel-reduced bases for the lattice $L \in \mathcal{L}_n$.

 Before we proceed, we observe that  the above result has a much easier variant. We can simply choose $\mathcal{B}$  at random from the 
 set $\mathfrak{S}_n$. In this model of random $\mathcal{B}$ the analogue of (*) is  quite straightforward to establish: it is part (iii) of the Theorem below. 
  However, it seems to us that the model of (*), i.e. first choosing a random lattice and then choosing
 a random basis for it,  is more natural (for example, in considering the behavior on LLL on many different
 bases for the same lattice). The  difficult part of this paper, then, is verifying that the two models are essentially equivalent,
 and this is accomplished by part (ii) of the Theorem below.

\begin{Theorem} \label{MainTheorem} 
As above, let  $N(L)$ be the number of Siegel-reduced bases for the lattice $L \in \mathcal{L}_n$, with reduction parameter $T$. 
Then:
\begin{itemize}
\item[(i)]  The $\mu_n$-expectation of $N(L)$ satisfies $\lim_n \frac{\log  \mathbb{E} N(L)}{n^3} = \frac{1}{6} \log T$.
\item[(ii)]  (Assuming the Riemann hypothesis):\footnote{We anticipate that this assumption should not be difficult to remove, but  it would require a more messy contour argument.}
The $\mu_n$-standard deviation of $N(L)$ is at most $\exp(-a n^2)$ times its mean. 

\item[(iii)] Fix $\delta > 0$ and let $X_{\delta} \subset \mathfrak{S}_n$
be the subset satisfying 
  $\|\mathbf{x}_n\|/ T^{(n-1)/2}< 1- \delta$. Then 
$$\frac{\mbox{measure} \ X_{\delta}}{\mbox{measure}  
\ \mathfrak{S}_n} \rightarrow 0$$
as $ n \rightarrow \infty$.  In words:
If we choose a basis from $\mathfrak{S}_n$ at random,  
the ratio $\|\mathbf{x}_n\|/ T^{(n-1)/2}$ is {\em concentrated at $1$.}  
\item[(iv)] Corollary to (i) and (ii):  in large enough dimension, $99.9\%$ of lattices have a  basis that is Siegel-reduced with parameter $1.0001$.

 \end{itemize} 
 \end{Theorem}

\proof
(i) is proved in  \S \ref{Siegelnumber} and (ii) is the main theorem (proved by the end of the paper). (iii) is proved in \S \ref{Siegelstatistics}.
  (iv) is an immediate consequence of (i) and (ii): for  any $T > 1$ (for example $T=1.0001$), the random variable $N(L)$  
 must be positive with the exception of a set of relative measure $\exp(-c n^2)$, for suitable $c > 0$. 
  
   \qed 
   
  While it is not at all surprising that $N(L)$ is concentrated  
 around its mean, the {\em extent} of the concentration is rather surprising.  To place the Gaussian $\exp(-an^2)$ that appears in (ii)
 in perspective,  let us note the following:  if we consider the set $\mathcal{L}' \subset \mathcal{L}_n$  of lattices  $L$ which possess a vector $\mathbf{y}$ of length $\|\mathbf{y}\| \leq 1$, 
 then in fact $\mu_n(\mathcal{L}') \geq e^{-C' n \log(n)}$ for suitable $C'$.  Now lattices in $\mathcal{L}'$ seem very atypical, because they possess a vector of very short length; one might expect this to strongly skew the set of Siegel-reduced bases  --
 but the theorem implies that they must mostly have the correct number of Siegel-reduced bases nonetheless! 
 
 We note finally that, in practice,
 LLL is often replaced by more sophisticated versions such as BKZ (see \cite{CN}, for instance).
 It would be interesting to try to understand the analogue of our results in that context.

 \subsection{Proof of the statement (*)} 
 
   Let us see how to derive the quoted statement on the previous page.

Fix $\delta > 0$.  For  $\mathcal{B}$ a Siegel-reduced basis, say that $\mathcal{B}$ is ``$\delta$-good" if $\|\mathbf{x}_n\|/T^{(n-1)/2} < 1-\delta$. (The larger the value of $\delta$, the better the basis!) 
 For each lattice $L$, let $N_{\delta}(L)$ be the  number of $\delta$-good Siegel-reduced bases for $L$ with parameter $T$. 
For $\alpha \in (0,1)$ consider now the exceptional set $\mathcal{L}_{!}$ of $L$ for which $N_{\delta}(L) > \alpha N(L)$. We will show its measure approaches zero.
For definiteness, let us show this measure is less than $\frac{1}{25}$. 

By (iii), for all sufficiently large $n$, we have
\begin{eqnarray} \int_{\mathcal{L}_n} N_{\delta}(L)  & =& \mbox{   measure of $X_{\delta}$} \\ &\leq& (\alpha/100) \int_{\mathcal{L}_n} N(L) \end{eqnarray}
For such $n$, then, we have
 $$   \int_{\mathcal{L}_{!}}N(L) < \frac{1}{100} \int_{\mathcal{L}_n} N(L).$$

 Because of  part (ii) of the Theorem,  there exists a subset $Y \subset \mathcal{L}_n$ 
 of measure $\leq \exp(-a n^2)$ so that 
 $N(L) >  0.5 \int_{\mathcal{L}_n} N(L)$ for $L \notin Y$.  Therefore,
\begin{eqnarray*} \int_{\mathcal{L}_{!} \backslash Y} N(L) < \frac{1}{100}  \int_{\mathcal{L}_n} N(L) &  \implies&  \mathrm{meas}(\mathcal{L_{!}} \backslash Y)  < \frac{1}{50}  \\
& \implies&  \mathrm{meas}(\mathcal{L}_{!}) < \frac{1}{50}+\exp(-a n^2) < \frac{1}{25}, \end{eqnarray*} 
where the last inequality again holds for large enough $n$. 
 \qed

It is quite easy to carry out the analog of (i) and (iii)   for LLL-reduction; it shows, for example, that the mean number $N_n$ of LLL-reduced
bases still satisfies $\lim \frac{\log N_n}{ n^3 \log T}= \frac{1}{6}$.  (For details, see \cite{KimThesis}.)  We conjecture that the other results
also remain valid for LLL reduced-bases.

 \subsection{About the proof} 
 
 The function $L \mapsto N(L)$ is, in the standard terminology of automorphic forms, a ``pseudo--Eisenstein'' series. 
 The theory of Eisenstein series allows one (at least in principle) to write it as an integral of standard Eisenstein series, 
 and then evaluate its $L^2$ norm.  So one ``just'' has to write everything out and bound each term. 
 
For $n=2$, this computation  amounts to computing the variance
in the number of vectors that a lattice has in a fixed ball. Such a computation
was apparently first done by W. M. Schmidt \cite{Schmidt} and rediscovered more recently by Athreya--Margulis \cite{AM}. 
But most of our complexity comes from the issues of large dimension. 
 
 The complication here is there are many types of ``standard Eisenstein series'' for $\SL_n$;
 they are indexed roughly by partitions of $n$. Correspondingly, the actual formula for the $L^2$-norm is very complicated.
 It involves many terms of the following type (cf. \eqref{prod2}):
 \begin{equation} \label{prod2} \prod_{b \in B}  \frac{   \xi(  \pm z + m_b) }{\xi(z+m_b +j_b)}\end{equation} 
for a certain purely imaginary variable $z$ and where $m_b, j_b$ are half-integers.  Also, $\xi$ is the completed Riemann zeta function,
$\xi(s) = \pi^{-s/2} \Gamma(s/2) \zeta(s)$. 

The only real concern  is that one might have $j_b < 0$ for many $B$.  In that case, the $\Gamma$-functions 
 in the numerator of \eqref{prod2}
are evaluated much further to the right than the $\Gamma$-functions in the denominator, and \eqref{prod2} would be ``large.''
This unfortunate possibility is ruled out by an explicit   combinatorial Lemma (Lemma \ref{CoreCombinatoriaLemma}),
which shows that in fact all the $j_b \geq 0$ (this is not at all obvious from the general presentation of the constant term of Eisenstein series,
although perhaps it is forced in some more subtle way by the internal structure of Eisenstein series). 

Besides this point, the other issues are minor. One needs plenty of careful book-keeping to keep track of measures on everything. 
 Another minor issue arises from the pole of the Riemann $\zeta$-function,  which we avoid by shifting to avoid it and using Cauchy's integral formula. 
 The use of the Riemann hypothesis arises in this step; it could likely be avoided with a little more care. 
 
It may be helpful to remember that, throughout the paper, terms of size $n^n$ will be essentially negligible.
Thus, for example, the total number of partitions of $n$ is negligible, compared to other quantities that we have to bound. 
We need to worry only about terms that are exponential in $n^2$ and higher.

In conclusion, from the point of view of automorphic forms, our result is in some sense
a  straighforward exercise. However, it seems to us that  the study of analysis on $\SL_n(\Z) \backslash \SL_n(\R)$ as ``$n \rightarrow \infty$''
is an interesting direction, and this paper represents  a first step in that direction.

\subsection{Acknowledgements}
Both authors acknowledge support from NSF grant   DMS-1401622. The first author (A.V.) was supported by a grant from the Packard Foundation.
The authors would like to thank Henry Cohn, Haseo Ki, and Stephen Miller for helpful discussions. 

 \tableofcontents

\section{Mathematical formulation of the main theorems}

Write $G=\SL_n(\R), \Gamma = \SL_n(\Z)$. 
 Note that 
 $ g \mapsto \Z^n g$
defines a homeomorphism between $\Gamma \backslash G$ and $\mathcal{L}_n$;
moreover, for fixed $g \in G$, 
the rule $$\gamma \mapsto \mbox{ the rows of $\gamma g$} $$
defines a bijection between $\Gamma$ and bases for the lattice $\Z^n g$.

Let $N,A,K$ be, as usual,
the subgroups of $G$ consisting of upper-triangular unipotent matrices,  
diagonal matrices with positive entries, and orthogonal matrices, respectively.    Let $B = AN $ be the Borel subgroup of upper triangular matrices. We write 
  $\Gamma_N = N \cap \Gamma$, etc.

 The product map gives a diffeomorphism
\begin{equation}  \label{iwa} N \times A \times K \stackrel{\sim}{\longrightarrow} G \end{equation}
For $g \in G$ we denote by $H(g)$ the ``$A$'' component, i.e. $g = n_g H(g) k_g$ with $n_g \in N, H(g) \in A, k_g \in K$.
Let  $\alpha_i : A \rightarrow \R_{+}$ is the simple root 
which sends $a =  \mathrm{diag}(a_1, \dots, a_n) \mapsto a_i/a_{i+1}$.
 (Note that $\prod a_i = 1$.)

Let $\mathbf{x}_1, \dots, \mathbf{x}_n$ for the rows of $g$, and 
  $\mathbf{x}_1^*, \dots, \mathbf{x}_n^*$  be the rows of $n_g^{-1} g = H(g) k_g$;
th $\mathbf{x}_i^*$ are orthogonal, $\| \mathbf{x}_i^*\| = a_i$,  and we have $\mathbf{x}_i = \mathbf{x}_i^* + \sum_{j > i} n_{ij} \mathbf{x}_j^*$. 
In particular, $\mathbf{x}_i^*$ must be the projection of $\mathbf{x}_i$ onto the orthogonal complement of $\langle \mathbf{x}_{i+1} , \dots, \mathbf{x}_n \rangle$.

We deduce that  the basis for $\Z^n g$ given by the rows of $g$ is Siegel-reduced with parameter $T$ if and only if
$$ \mbox{$\alpha_i(H(g)) \geq T^{-1}$, and for $i \neq j$ we have  $|(n_g)_{ij}| \leq 1/2$},$$

Let $f$ be the function on $A$ given by
\begin{equation} \label{fdef} f(a) = \prod  1_{\alpha_i > 1/T} \end{equation} 
where $1_{\dots}$ denotes ``characteristic function''.

The pseudo-Eisenstein series induced from $f$ is by definition
\begin{equation} \label{pseudo-Eis-def} E_f(g) = \sum_{\gamma \in  \Gamma_B \backslash \Gamma} f(H(\gamma g))\end{equation}
 
Now $g \mapsto \Z^n g$ yields (away from a measure zero set of $g$) 
a bijection between $\Gamma_N \backslash \Gamma$
and  reduced bases of $\Z^n g$.
Since $[\Gamma_B: \Gamma_N] = 2^{n-1}$, we see that 
$$ E_f(g) = \frac{1}{2^{n-1}} N(\Z^n g),$$
where $N(\Z^n g)$ was as before the number of  Siegel-reduced bases  with parameter $T$ inside the lattice $\Z^n g$.
(Note that $f$ implicitly depends on $T$). 
%

  Let $\overline{E_f}$
be the average value of $E_f$ over the space of lattices. We will prove the following (with $\xi(s)$ the completed Riemann zeta function,
$\xi(s) = \pi^{-s/2} \Gamma(s/2) \zeta(s)$): 

\begin{equation} \label{a}  \overline{E_f} = T^{\frac{1}{6}(n^3-n)} \cdot  \frac{1}{n} \prod_{i=2}^n \frac{1} {i(n-i)}  \cdot  \xi(2)^{-1} \dots \xi(n)^{-1}  \end{equation} 
\begin{equation} \label{c} \frac{ \| E_f - \overline{E_f} \|_2  }{ \|E_f\|_2} \leq A e^{-\delta n^2}\end{equation} 
for suitable constants $A, \delta$.   Note that $A$ will depend on $T$, but $\delta$ does not.  Said differently,  the orthogonal projection of $E_f$ to the orthogonal complement of constants
accounts for  at most $ A e^{-\delta n^2}$ of the $L^2$-norm. 

The proof of \eqref{a} is quite straightforward and is completed in 
\S \ref{Siegelnumber}, but \eqref{c} is quite a bit deeper. It uses the full spectral theory of automorphic forms on $\SL_n$. 

\section{Setup}

%
%
%

\subsection{Haar measures} 
 
Fix Haar measure $dn$ on $N$ such that the covolume of $\Gamma_N$ is $1$. 
 Explicitly we take $dn = \prod_{i < j} dn_{ij}$. 
 We equip $A$ with the Haar measure $da := \prod d \alpha_i/\alpha_i$.   

Let $2\rho: A \rightarrow \mathbf{R}_{+}$ be the sum of all positive roots of $A$. We will often use additive, rather than multiplicative,
notation for characters of $A$; therefore,  
\begin{equation} \label{rhodef} 2 \rho =  \sum_i  {i(n-i)} \alpha_i.  \end{equation}

By means of $NA \simeq G/K$ we get a $G$-invariant measure    on $G/K$: \begin{equation} \label{Gmeasure} dn \cdot da \cdot a^{-2\rho}. \end{equation} 
We pull this measure back to $G$ via $G \rightarrow G/K$, normalizing the measure of $K$ to equal $1$. 

\subsection{The vector subspace $\sum x_i = 0$} \label{3meas}
Let $U_n = \{\sum x_i = 0 \} \subset \mathbf{R}^n$. This subspace
has three natural measures on it;

\begin{itemize}
\item The ``fibral'' measure $\nu_f$ given by disintegrating Lebesgue measure over the map
$\mathbf{R}^n \rightarrow \mathbf{R}$ given by $(x_i) \mapsto \sum x_i$. 
This is given as a differential form by $|\omega|$, where  $\omega = \wedge_{i=1}^{n-1} dx_i$ (or indeed the same product omitting any one of $dx_1, \wedge, \dots, dx_n$). 

 \item The ``Riemannian'' measure $\nu_R$, corresponding to the restriction of the standard inner product on $\mathbf{R}^n$. 
 
\item The ``big'' measure $\nu_b$, given by $|\omega'|$ with $\omega' = \wedge_{i=1}^{n-1} dx_i-dx_{i+1}$, i.e. 
Lebesgue measure if we identify $U \simeq \mathbf{R}^{n-1}$ via $(x_2-x_1, \dots, x_{n}-x_{n-1})$. 

\end{itemize}
For example, when $n=2$, the measure of the set $\{ x_1+x_2=0, 0 \leq x_1 \leq 1\}$
is $1, \sqrt{2}, 2$ according to the three measures in the order specified above.

These measures are related by
$$\nu_b = \sqrt{n} \nu_R,  \nu_f = \frac{1}{\sqrt{n}} \nu_R$$

Identify $U_n$ with its own dual via $\langle x_i, y_i \rangle = \sum x_i y_i$. 
Then (for $dx$ any of the measures just noted) we can define the Fourier transform of a function $f$ on $U_n$ via
$$\hat{f}(k) = \int f(x) e^{2 \pi i \langle x, k \rangle} dx, \ \ (y \in U_n) $$ 
and there is a Fourier inversion formula, replacing  then $dx$ by the dual measure. 
The dual measure to $\nu_b$ is $\nu_f$, and vice versa; the measure $\nu_R$ is self-dual.

\subsection{Volume}

The volume of $\Gamma \backslash G$
in our normalization  is  \footnote{Note if we denote by $\mu_T$ the ``Tamagawa'' measure, the measure induced by integrally normalized differential form,   it is related to our measure $\mu$ via
$ \mu_T = \frac{1}{n} \prod_{j=2}^{n}  S_j  \cdot \mu $
where $S_j = \frac{2 \pi^{j/2}}{\Gamma(j/2)}$.  Here we used
the fact that the measure of $\mathrm{SO}_n$ is the product $S_2 \dots S_n$.
 }
\begin{equation} \label{mentioned} \vol(\Gamma \backslash G) =n \prod_{j=2}^{n}  \zeta(j) S_j^{-1} = \frac{n}{2^{n-1}}  \underbrace{\xi(2) \dots \xi(n)}_{Q_n},\end{equation}  where
$\xi$ is the completed $\zeta$-function and $S_j = \frac{2 \pi^{j/2}}{\Gamma(j/2)}$ is   the surface area of the sphere in $\mathbf{R}^j$.  We write for short 
\begin{equation} \label{Qdef} Q_n = \xi(2) \dots \xi(n)\end{equation} 

\subsection{Characters of $A$ and their parameterization} \label{chiPar}
   
    Let $\chi$ be a character of $A$, i.e. a character $A \rightarrow \mathbf{C}^*$. 
    
    We will usually represent $\chi$ in one of two ways:
    $$ \chi(a) = \prod_{i=1}^n a_i^{\nu_i}, \ \ \sum \nu_i= 0,$$
    or  as $$\chi = \sum_{i=1}^{n-1} \mu_i \alpha_i, \mbox{ i.e. $\chi(\diag(a_1, \dots, a_n))= \prod_{i=1}^{n-1} (a_{i}/a_{i+1})^{\mu_i}$.} $$
    
    We will often write $\nu_{ij}$ as a shorthand for $\nu_i-\nu_j$.
    Also we will often write $\nu_i(\chi)$ or $\mu_i(\chi)$ for the parameters $\mu_i$ or $\nu_i$ as above; they are related via 
    $$\nu_1=\mu_1, \nu_2=\mu_2-\mu_1, \dots, \nu_{n-1} = \mu_{n-1}-\mu_{n-2}, \nu_n= - \mu_{n-1}$$
    and in the reverse direction
\begin{eqnarray} \label{reverse_formula} \mu_{n-1} &=& - \nu_n = \nu_1 + \dots \nu_{n-1}  , \\ 
\nonumber  \mu_{n-2} &=&  - \nu_{n-1} - \nu_{n} = \nu_1 + \dots + \nu_{n-2}, \dots
\end{eqnarray}
    and so on. Note that 
     \begin{equation} \label{eqmes}  d\nu_1 \wedge \dots \wedge d\nu_{n-1} = \pm  d\mu_1 \wedge \dots \wedge d\mu_{n-1}.\end{equation}
     
   We will write \begin{equation} \label{chinotn} \wt(\chi) := \sum \mu_i(\chi), \mathrm{P}(\chi) = \prod \mu_i(\chi). \end{equation}
 e.g. from \eqref{rhodef} we have $\wt(2\rho) = \sum_{i=1}^{n-1} i(n-i) =\frac{1}{6}(n^3-n)$.

     
From \eqref{fdef} we get, for any character $\chi = \prod \alpha_i^{\mu_i}$: 
\begin{equation} \label{frofro}  \int_{A} f \chi^{-1}  d a = \prod_{i=1}^{n-1} \int_{1/T}^{\infty} u^{-\mu_i} \frac{du}{u} =  \frac{1}{\prod \mu_i} T^{\sum \mu_i} =  \frac{T^{\wt(\chi)}}{\mathrm{P}(\chi)} \end{equation}
(this computation is the main reason  for introducing the notation \eqref{chinotn}).

%
%
%
%
%
%
%
%
%
    
\subsection{The average number of Siegel reduced bases} \label{Siegelnumber}

From our previous discussion it follows that the average   $\overline{E_f} $ of  $E_f$ is given by $\frac{1}{\vol(\Gamma \backslash G)} \int E_f$.
 By unfolding, $\int E_f =  2^{1-n} \int a^{-2\rho} f(a) $:
\begin{equation} \label{unfold} \int E_f =  \int_{\Gamma \backslash G}  dg  \left( \sum_{\Gamma_B \backslash \Gamma} f(H(\gamma g) \right)  = \int_{\Gamma_B \backslash G} f(H(g))  dg   
\end{equation}
$$ = 2^{1-n} \int_{\Gamma_N \backslash G} f(g) dg \stackrel{\eqref{Gmeasure}}{=} 
 2^{1-n} \int_{} f(a) a^{-2\rho} da $$
  We conclude from \eqref{frofro} and \eqref{rhodef}: 

$$ \overline{E_f} = T^{\frac{1}{6}(n^3-n)} \cdot  \frac{1}{n} \prod_{i=1}^{n-1} \frac{1} {i(n-i)}  \cdot  \xi(2)^{-1} \dots \xi(n)^{-1}.$$

\subsection{The mean length of the first vector of a Siegel-reduced basis} \label{Siegelstatistics} 

By a similar computation we can compute the mean value of $a_n$:
Note $$ \alpha_1 \alpha_2^2 \dots \alpha_{n-2}^{n-2} \alpha_{n-1}^{n-1} = (a_1 \dots a_{n-1})/a_n^{n-1} = a_n^{-n}.$$
So we want to compute the mean value of $\prod_{i=1}^{n-1} \alpha_i^{-i/n}$. We do this just as in \S \ref{Siegelnumber}; the mean value equals: 
$$ \frac{T^{\sum i(n-i) + i/n}}{T^{\sum i(n-i)}} \cdot \prod \frac{i (n-i)+i/n}{i(n-i)} \sim T^{\frac{n-1}{2}}.$$
where $\sim$ here means that the ratio of both sides approaches $1$ as $n \rightarrow \infty$. 
Note that this mean value corresponds exactly to the ``worst case behavior'' of $\|\mathbf{x}_n\|$ from \eqref{worst-case}.

In fact $\log(a_n)$ is obtained  as the convolution $\sum \frac{i}{n} Y_i$, 
where each $Y_i$ is a variable on $[\log(T), \infty]$
with distribution function proportional to $e^{-i (n-i) Y}$. It easily follows that in fact $\log(a_n)$ is concentrated around its mean value.

%
%
%
%

    \subsection{Levi subgroups and their parameterization} \label{levis}
    Given $n_1, \dots, n_k$ such that  $n =\sum_{i=1}^{k} n_i$ 
    let  $M$ be the corresponding Levi  subgroup of block diagonal matrices, with Levi subgroup $\mathrm{S} (\prod_{i=1}^k \GL_{n_i}) \subset \SL_n$;
    we write
    $$N_M = N \cap M, K_M = K \cap M.$$
 We  also  write $$N_1=n_1, N_2 =n_1+n_2, N_3=n_1+n_2+n_3, \dots$$
 
 The {\em blocks} of the Levi subgroup (thinking of it as block diagonal matrices) are parameterized by the intervals
\begin{equation} \label{blockdef} [1, N_1], \ \ [N_1+1, N_2], \ \ , [N_2+1, N_3], \dots \end{equation}

We call  the tuple $n_1, \dots, n_k$, or equivalently the  standard Levi subgroup $M$, {\em good} if $n_1 \leq n_2 \leq \dots$. Among each associate class of Levi subgroups there is a good representative; for this reason,
it will be enough for us to consider only good Levi subgroups. 

We have a decomposition
\begin{equation}  \label{factor} M = M_0 Z_M\end{equation} 
where $M_0 \subset M $ is the subgroup of elements $\mathrm{S} (\prod_{i=1}^k \GL_{\pm 1})$  with determinant $\pm 1$ on each factor; and $Z_M$
is the subgroup of matrices that are positive scalar in each block (so $Z_M \simeq \mathbf{R}_{>0}^{k-1}$); a typical element of $Z_M$ looks like
$$ \diag \left( \underbrace{z_1, z_1, \dots, z_1}_{n_1}, \underbrace{z_2, \dots, z_2}_{n_2}, \dots , z_k \right)$$

 We equip $Z_M$ with the Haar measure $\prod d\beta_i/\beta_i$, where $\beta_i =  \alpha_{N_i}$ (considering $Z_M \subset A$).    Inside $M_0$ is the torus $A_M:= A \cap M_0$, which we equip with the analog of the measure $da$ on $A$, namely the measure
 $$ \prod_{j} \frac{d \alpha_i}{\alpha_i}$$
 where we take the product only over those $i$ not equal to any of $N_1, N_2, \dots$. 
With these normalizations, the product map $Z_M \times A_M \rightarrow A$ preserves measures. 

\subsection{More on measures}
 
Let $P_M$ be the parabolic subgroup geneated by $B$ and $M$. 
Let $U_{P_M}$ (or just $U_M$ for short) be the unipotent radical of $P_M$; thus, 
$P_M = U_M \cdot M$, and $U_M \cdot N_M = N$, where $N_M = N \cap M$. 

We normalize measure on $N_M$
so that $\vol(N_M/N_M \cap \Gamma) = 1$.  We normalize measure on $U_M$ in exactly the same way.  

This induces a normalization of measure on $M/K_M$, via
$$U_M \times M/K_M \simeq G/K$$
-- recall that the measure on $G/K$ came from $NA$, see \eqref{Gmeasure}.
More explicitly, this is the measure obtained via
$$M  = N_M A  K_M$$
by equipping $N_M$ with the measure where $\vol(N_M/N_M \cap \Gamma)=1$,
and $A$ with the same measure as before. The factorization \eqref{factor},
induces
$$M/K_M \simeq M_0/K_{M} \times A_M$$
 
 \subsection{Characters of $M$}  \label{CharM}
 Continue with the notation for $M$ as previous; in particular, $k$ is the number of blocks. 
Let $\mathfrak{a}_M^*$ be the space 
$$ (\nu_1, \dots, \nu_k) \in \C^k: \sum n_i \nu_i = 0.$$ 
We say that $\nu$ is {\em unitary} if every $\nu_i$ belongs to $i \mathbb{R}$. 
We denote the subset of $\mathfrak{a}_M^*$ as $\mathfrak{a}_{M,0}^*$:
$$  \mathfrak{a}_{M,0}^* = \{ (\nu_1, \dots, \nu_k) \in (i \R)^k: \sum n_i \nu_i = 0\}.$$

Elements $\nu \in \mathfrak{a}_M^*$ parameterize characters $\chi_{\nu}$ of $M$ defined by 
\begin{equation} \label{chinudef} \chi_{\nu}:  (g_1 \in \GL_{n_1} , g_2 \in \GL_{n_2}, \dots) \mapsto \prod_{i=1}^k |\det g_i|^{\nu_i}\end{equation}
these are precisely those characters of $M$ that are trivial on $M \cap K$, which is all that is of interest for us. 

On $Z_M$ this character is given by 
\begin{equation} \label{Zmchar} \chi_{\nu}(z_1, z_1, \dots,  , z_k, z_k) = \prod z_i^{n_i \nu_i} \end{equation}
Since   $z_k = \prod_{i  < k } z_i^{-n_i/n_k}$ we  
 can rewrite this as
$$\chi_{\nu}(z_1, z_1, \dots,  , z_k, z_k)  = \prod_{i < k} z_i^{n_i (\nu_i-\nu_k)}$$

The measure on $Z_M$ is equal to $\frac{n}{n_k} \frac{dz_1}{z_1} \dots \frac{dz_{k-1}}{z_{k-1}}$.
We may identify $\mathfrak{a}_{M,0}^*$ with the dual group to $Z_M$, and 
we want to compute the corresponding dual measure. 

From  $\sum n_i \nu_i = 0$ we deduce that $\nu_k = \sum_{i=1}^{k-1} - \frac{n_i}{n_k} \nu_i$, and so  
$$d(\nu_1-\nu_k) \wedge \dots \wedge d(\nu_{k-1}-\nu_k) =  \left( 1 +\sum_{i=1}^{k-1} \frac{n_i}{n_k} \right) d\nu_1 \wedge \dots \wedge d\nu_{k-1}
= \frac{n}{n_k}  d\nu_1 \wedge \dots \wedge d\nu_{k-1}  $$
We deduce correspondingly that the measure on $\mathfrak{a}_{M,0}^*$ that is dual to $Z_M$ is given by (the absolute value of)
\begin{equation} \label{numeasure} \frac{n_k}{n} \prod_{i=1}^{k-1} \frac{ n_i  d(\nu_i -\nu_k)}{2 \pi i} = \prod_{i=1}^{k-1} \frac{ n_i  d\nu_i}{2 \pi i}. \end{equation}

 

 \subsection{Volume of quotients for Levi subgroups}
 We will need to compute the volume $V_M$  of $\Gamma_M \backslash M_0/K_M$. 
It equals  the product of the quantity \eqref{mentioned} over $n=n_1, \dots, n_k$. In particular, if we write  
\begin{equation} \label{Mvoldef}   Q_M :=   \prod  Q_{n_i}\end{equation}  
then $e^{-A n \log(n)} \leq \frac{V_M}{Q_M} \leq e^{A n \log(n)}$ for a suitable constant $A$.

%
%
%
%
%
%
%
%

\section{The principal Eisenstein series and its constant term}  

We will summarize what we need from the theory of Eisenstein series, presented in ``classical'' language. 
A clear summary of the theory of Eisenstein series, but in adelic language, is given in \cite{Arthur}: the computation of constant terms is Lemma 7,
and the holomorphicity of Eisenstein series on the unitary axis is stated in the {\em Main Theorem}.  A reference which uses 
the language of real groups and is closer to our presentation here is \cite{Langlands}. 
 
 \subsection{Borel Eisenstein series} 

 The usual Borel Eisenstein series is indexed  by characters $\chi $ of $A$:
 $$ E_B(\chi, g) = \sum_{\gamma \in \Gamma_B \backslash \Gamma} \langle \chi + \rho, H(\gamma g) \rangle$$
 The constant term of $E_B$ -- i.e., $ \int_{\Gamma \cap N\backslash N} E_B(ng) dn$  -- depends only on the ``$A$'' component of the $NAK$ decomposition, and is of the form 
\begin{equation} \label{constant term} (E_B)_N =  \sum_{w \in W} a^{w \nu + \rho}  \prod_{w \alpha < 0, \alpha > 0} \frac{\xi(s_{\alpha})}{\xi(s_{\alpha}+1)}\end{equation} 

 where $W \simeq S_n$ is the Weyl group acting by coordinate permutation on $A$, $\alpha$ ranges over coroots,  $s_{\alpha} = \langle \chi, \alpha^{\vee} \rangle$ and $\xi(s) = \pi^{-s/2}\Gamma(s/2) \zeta(s)$  is the completed $\zeta$-function; more explicitly, in the coordinates $\chi = (\nu_1, \dots, \nu_n)$ introduced in \S \ref{chiPar}, we have 
 \begin{equation} \label{constant term2} (E_B)_N =  \sum_{\sigma \in S_n} a^{\sigma \nu + \rho}  \prod_{i<j, \sigma(i) > \sigma(j)} \frac{\xi(\nu_{ij})}{\xi(\nu_{ij}+1)}\end{equation} 
where we wrote $\nu_{ij} = \nu_i - \nu_j$, and also   $\sigma \nu = \left(\nu_{\sigma^{-1}(1)}, \nu_{\sigma^{-1}(2)}, \dots  \right)$ or more evocatively $(\sigma \nu)_{\sigma(i)} = \nu_i$
.

 In this normalization, the ``unitary axis'' is given by $\mathrm{Re}(\chi) = 0$, i.e. $\mathrm{Re}(\nu_i) = 0$.   The point $\chi =  \rho$ is the intersection of all lines $s_{\alpha} =1$
 for all {\em simple} roots $\alpha$, i.e. the intersection of lines $\nu_i - \nu_{i+1}=1$.  If we take iterated residues of $E$ along all of these hyperplanes,  the only term that contributes is $\sigma \in S_n$
 given by $i \mapsto n+1-i$; we get 
\begin{equation} \label{resEcomp}  \mbox{residue of $E$} = \mbox{ constant function with value $Q_n^{-1}$.}\end{equation}
%
%
 \subsection{Degenerate Eisenstein series} \label{sec:degenerate}
 
 We will be also interested in the Eisenstein series induced from  one-dimensional representations of a  Levi subgroup.

 Let $M$ be a standard Levi factor (cf. \S \ref{levis})  corresponding to the  decomposition $ \sum_{i=1}^r n_i = n$; 
 we will identify this with the partition of \begin{equation} \label{partitionM} \{1,\dots ,n\}
 = [1,n_1] \cup [n_1+1, n_1+n_2] \cup \dots [\sum_{i=1}^{r-1} n_i + 1, n]. \end{equation}  
Just to recall our notation, we will refer to the subsets appearing above -- $[1,n_1]$, $[n_1+1, n_1+n_2]$ and so forth -- as the {\em blocks} associated with the Levi $M$;
and we put $N_1=n_1, N_2=n_1+n_2$, etc. 
 
 We will only consider {\em good} Levi subgroups with $n_1 \leq n_2 \leq \dots$. 
  Let $P =MB$ be the corresponding parabolic and $\Gamma_P = \Gamma \cap P$. 
 
 Let $\nu \in \mathfrak{a}_M^*$ be as in \S \ref{CharM}.
 We define the degenerate Eisenstein series $E_{M, \nu}$ parameterized by $\nu$, by means of
  $$E_{M,\nu}(g) = \sum_{\gamma \in \Gamma_P \backslash \Gamma}  \langle \chi_{\nu} + \rho_P, H(\gamma g) \rangle$$
 where $\chi_{\nu}$ is as in \eqref{chinudef}, $ \rho_P$  is the character of $M$ defined by
 $$ \rho_P:  (g_1 \in \GL_{n_1} , g_2 \in \GL_{n_2}, \dots) \mapsto   \prod_{i< j} \left( \frac{ |\det g_i|^{n_j} }{ |\det g_j|^{n_i}}\right)^{1/2}.$$
 
\subsection{The constant term of the degenerate Eisenstein series $E_{M, \nu}$}
We now want to compute the constant term of $E_{M,\nu}$. We do this by interpreting it as a residue of $E_{B}(\chi)$. 
For this we use ``induction in stages,''  
  we can express $E_B(\chi)$ as the iterated Eisenstein series $E_P^G E^{M}_{M \cap B}(\chi)$.
  
More precisely:
$$ \sum_{\gamma \in \Gamma_B \backslash \Gamma} \langle \chi+\rho, H(\gamma g)  \rangle = \sum_{\gamma_2 \in \Gamma_P \backslash \Gamma} \sum_{\gamma_1 \in \Gamma_B \backslash \Gamma_P}
\langle \chi + \rho, H(\gamma_1 \gamma_2 g) \rangle
 $$ and we have an identification $\Gamma_{M \cap B} \backslash \Gamma_M \simeq \Gamma_B \backslash \Gamma_P$.
The function $$g \mapsto \sum_{\gamma_2 \in \Gamma_B \backslash \Gamma_P} \langle \chi+\rho, H(\gamma_1 g) \rangle,$$
defines a function on $\Gamma_P U_P \backslash G/K \simeq \Gamma_M \backslash M /K_M$, which we call $E^M_{M \cap B}(\chi)$;
this is an Eisenstein series for $M$. 
By an analogue of \eqref{resEcomp}, the iterated residue  of $E^M_{M \cap B}(\chi)$  along the lines 
\begin{equation} \label{tier} \nu_i-\nu_{i+1}=1 , \ \ \mbox{ for all } i \notin \{N_1, N_2 , \dots \}\end{equation}  gives the constant function with value $Q_M^{-1}$  
(see \eqref{Mvoldef}), and therefore
 $$E_{M,\nu}  =  Q_M \mathrm{Res} \   E_B(\upsilon)$$
where  the residue is still taken along \eqref{tier}, and we defined the shifted parameter:   {\small  \begin{equation} \label{Mchar}  \upsilon = \left( \frac{n_1-1}{2}+ \nu_1, \frac{n_1-3}{2} +\nu_1, \dots, -\frac{n_1-1}{2} + \nu_1 , \frac{n_2-1}{2} +\nu_2, \dots, - \frac{n_2-1}{2} + \nu_2, \dots \right). \end{equation} }

By means of  \eqref{constant term} we get the following expression for the constant term of $E_{M, \nu}$:
Define
\begin{equation}
\label{SMdef} S[M] = \{ \sigma \in S_n: \mbox{$\sigma$  
is monotone decreasing on each of $[1, n_1]$, $[n_1+1, n_2]$ etc.}\}\end{equation}
  Therefore
the term $ \prod_{i<j, \sigma(i) > \sigma(j)} \frac{\xi(\nu_{ij})}{\xi(\nu_{ij}+1)}$ that appears in \eqref{constant term2}  has nonzero iterated residue along \eqref{tier} precisely when $\sigma \in S[M]$, because to have nonzero residue along $\nu_{i,i+1} = 1$
we must contain a term corresponding to $(i, j=i+1)$, i.e. we must have  $\sigma(i) > \sigma(i+1)$, i.e. $\sigma$ must be decreasing on each block $[1,n_1], [n_1+1, n_2]$ and so on. 
Therefore, 

\begin{eqnarray} \nonumber  \left( E_{M,\nu} \right)_N &=& Q_M \sum_{\sigma \in S[M]}  \mathrm{Res}_{\eqref{tier}} \prod_{i < j, \sigma(i) > \sigma(j)}    \frac{\xi(\upsilon_{ij})}{\xi(\upsilon_{ij}+1)}  a^{\sigma \upsilon  + \rho}
\\  \label{Cterm}  &=& \sum_{\sigma \in S[M]} \prod_{i < j, \sigma(i) > \sigma(j), i \nsim j}   \frac{\xi(\upsilon_{ij})}{\xi(\upsilon_{ij}+1)}  a^{\sigma \upsilon  + \rho}
\end{eqnarray}  where  as before $\upsilon_{ij} = \upsilon_i-\upsilon_j$ and $\upsilon$ is as in \eqref{Mchar}; also $i \nsim j$ means they are different parts of the partition \eqref{partitionM} corresponding to $M$.
For the second line, we used the fact that the contribution of all terms  with $(i,j)$ in the same block is precisely $Q_M^{-1}$.

In summary, the constant term of the degenerate Eisenstein series $E_{M,\nu}$ is given by \eqref{Cterm}, where $\upsilon$ is given in terms of $\nu$ in \eqref{Mchar}; this expression is very similar
to the constant term \eqref{constant term2} of the full Eisenstein series $E_B(\chi)$.

{\em Important warning.} 
  In general, the terms in this expression can still have poles; a priori, they determine only a meromorphic function of  $\nu$,   and \eqref{Cterm} is valid as an equality of meromorphic functions of $\nu$. However, these meromorphic functions are 
 necessarily holomorphic on the line $\mathrm{Re}(\nu_{i}) = 0$: a basic result of the theory of Eisenstein series  (see e.g. the first sentence of the Main Theorem, \cite{Arthur}) is that the Eisenstein series
 induced from a discrete-series representation is holomorphic on the ``critical line,'' which in this case corresponds to $\mathrm{Re}(\upsilon_{ij}) = 0$.

 It is more convenient to rewrite this in a way that is indexed by {\em blocks} of the Levi $M$. Recall the blocks are just the intervals
 of integers corresponding to blocks of the Levi, cf. \eqref{blockdef}. 
We can write: 

\begin{equation} \label{constant2}  \left( E_{M,\nu} \right)_N =  \sum_{\sigma \in S[M]}  \prod_{A  <  B} \prod_{i \in A, j \in B, \sigma(i) > \sigma(j) }    \frac{\xi(\upsilon_{ij})}{\xi(\upsilon_{ij}+1)}  a^{\sigma \upsilon  + \rho}
 \end{equation}
 Here $A< B$ means that $A$ precedes $B$ in the natural ordering. 
\subsection{Spectral theory} 
 Let $M$ be a standard Levi subgroup, as above.  
 
 Let $W_M$ be the group of self-equivalences of $M$, that is to say, the set of $w \in W$ 
 with the property that $w$ preserves the center of $M$.
   For example, if  the lengths of blocks $n_1, n_2, \dots$ are all pairwise disjoint, then $W_M$ consists of those elements $w \in S_n$ 
 which stabilize, setwise, each block $[1,n_1]$, $[n_1+1, n_1+n_2]$ and so forth.
%
 
 Let $f$ be a function on $A$, and form the pseudo-Eisenstein series $E_f$ as in \eqref{pseudo-Eis-def}.
   We have a spectral decomposition 
\begin{equation} \label{SpectralDecomposition} E_f = \sum_{M} V_M^{-1} \frac{1}{|W_M|} \int_{\nu \in \mathfrak{a}_{M,0}^*} d\nu   \langle E_f,  E_{M,\nu} \rangle  E_{M,\nu}\end{equation} 
   where the sum is taken over  good Levi subgroups $M$ (see \S \ref{levis}), the group $W_M$ is as above,  the $\nu$-integral is taken over\ $\mathfrak{a}_{M,0}^*$,    
and the measure to be taken on the space of parameters is that dual to the Haar measure on $Z_M$  (cf. \S \ref{levis} and \eqref{numeasure}). 
   
 Let us explain briefly the origin of \eqref{SpectralDecomposition}. Indeed, there is a corresponding expansion for any function $\varphi$ on $\Gamma \backslash G$, but it in general involves
Eisenstein series $E_{M, \nu}(\psi)$  induced from all  automorphic forms $\psi$ on $M$ that lie in $L^2$ modulo center.   We must check that only those
$\psi$ that arise from characters of $M$ contribute to this spectral expansion (and then only those characters trivial on $M \cap K$ give a nonzero contribution;
these are precisely the characters of
\eqref{chinudef}). 
In other words, we must verify that $\langle E_f, E_{M,\nu}(\psi) \rangle$ vanish unless $\psi$ is a character of $M$. 
By unfolding, this inner product vanishes unless $E_{M,\nu}(\psi)$ has nontrivial constant term along $N$,
which implies in particular that $\psi$ itself has nontrivial constant term along $N_M$. 
 So it is enough to verify that, on $\GL_b$, any automorphic form
 of the discrete series, with nontrivial constant term along the unipotent radical of a Borel subgroup,
 must in fact be a character. This is not a triviality; it follows from the computation of the discrete spectrum of $\GL$ by Moeglin and Waldspurger \cite{MW};
 they show that the discrete spectrum for $\GL_n$ arises from  a divisor $a|n$ and a cusp form $\pi_a$ on $\GL_a$:
 one takes a certain residual Eisenstein series  $\Pi_{\mathcal{E}}$ induced from $\pi_a \boxtimes \pi_a \cdots \boxtimes \pi_a$.  In particular,
 if $a \neq 1$, the constant term along the unipotent radical of a Borel subgroup involves a constant term for $\pi_a$, and is zero;
 in the case $a=1$, then $\Pi_{\mathcal{E}}$ is one-dimensional.

   Note that the inner product $\langle E_f, E_{M,\nu} \rangle_{\Gamma \backslash G}$
   can be computed by unfolding $E_f$: it equals 
   \begin{equation} \label{unfold2} \langle E_f, E_{M,\nu} \rangle_{\Gamma \backslash G} = 2^{1-n}  \langle  f ,   a^{-2\rho} (E_{M,\nu})_N \rangle_A\end{equation}
   by the same argument as \eqref{unfold}; the inner product is computed in $L^2(A)$.

  \subsection{Re-indexing} \label{sec:reindexing}

  There is a bijection
  $$\mbox{ divisions $\mathcal{J}$  of 
  $\{1, \dots, n\}$} \longrightarrow  \mbox{  Levi subgroups $M$, $\sigma \in S[M]$ }  $$
 where by a ``division'' $\mathcal{J}$ what we mean is an ordered collection of disjoint 
  subsets   \ $J_1, \dots, J_k \subset \{1, 2, \dots, n\}$,   where $\coprod J_i = \{1, \dots, n\}$. 
  
  The bijection is defined thus:  associated to $\mathcal{J}$
  we take that Levi subgroup with $n_1 = \# J_1 , n_2 = \# J_2$ and so on. 
Write 
   $J_k = \{j_{k,1}, \dots, j_{k,n_k}\}$ with $j_1 < \dots < j_{n_k}$; 
   then there is a unique element  $\sigma_{\mathcal{J}}$ of $S[M]$ 
   where we take
   $$\sigma_{\mathcal{J}} :  \{N_{k}+1, \dots, N_{k+1}\} \rightarrow J_k$$
   but reversing order, that is to say:
   $$\sigma_{\mathcal{J}}(N_{k+1}-t)  = j_{k, t+1} \ \ \  0 \leq t \leq n_{k+1}-1.$$

 Note, for later use, 
that the number of such $\mathcal{J}$ as above is clearly at most $n^n$
  (clearly $k \leq n$ and the division is describeed by a function
  $\{1, \dots, n\} \rightarrow \{1, \dots, k\}$).  In particular,
  $$ \mbox{ number of possible $\mathcal{J}$} \times n! \leq n^{C n}$$
  for suitable $C$.

\subsection{The $\rho$ parameters} \label{rhoJdef}
Suppose $(M, \sigma \in S[M])$ corresponds to $\mathcal{J}$, as above. 
We make some computations related to the half-sum of positive roots for $M$. 

Let $\rho_M$  be the half-sum of positive roots for $M$, that is to say,  $2\rho_M = \sum e_{ij}$ 
 over all roots $e_{i,j}$ where $i,j$ belong to the same block of the partition defined by $M$, and also $i < j$. 
Define $\rho_{\mathcal{J}}$ by the rule
 $$2 \rho_{\mathcal{J}} = - \sigma (2 \rho_M)$$
 Therefore,  $$2\rho_{\mathcal{J}} = \sum_{i \sim_{\mathcal{J}} j, i < j } e_{ij} $$
where now $i \sim j$ means that they belong to the same part of the partition defined by $\mathcal{J}$.

%
%
Note that if $\mathcal{J}$ corresponds to $\sigma \in S[M]$, and $\upsilon$ is a normalized character corresponding
to {\em unitary}  $\nu$,  as in \eqref{Mchar}, then  $\mathrm{Re}(\upsilon) = \rho_M$ and thus 
\begin{equation}  \label{realpart} \mathrm{Re} \left(  \sigma \upsilon  \right) = - \rho_{\mathcal{J}}\end{equation}

 We will compute the $\mu_i$-coordinates of $\rho_{\mathcal{J}}$ for later use.
Visibly 
 $\mu_i(2 \rho_{\mathcal{J}})$ counts the number of pairs $a \sim b$ with $a \leq i$ and $b > i$. 
 Therefore,
 
$$ \wt(2\rho_{\mathcal{J}})  = \sum_{i \sim j, i <j} (j-i)$$  
$$ \wt(2\rho-2\rho_{\mathcal{J}}) = \sum_{i \nsim j, i< j } (j-i)$$
 
 Note that, if $\mathcal{J}$ has more than one part,  then
\begin{equation} \label{muicomb} \mu_i(2\rho-2\rho_{\mathcal{J}}) \geq 1 \mbox{ for all $i$.} \end{equation} 
 Indeed, if $\mu_i(\rho-\rho_{\mathcal{J}}) $ were zero, it means that every $a \leq i$ and $b> i$ belong to the same part of $\mathcal{J}$,
 which forces that $\mathcal{J}$ has just one part.

\begin{Lemma} \label{weightbound}
  Let $\noJ$ be the number of parts of $\mathcal{J}$
and let $z$ be the size of the largest part. For all $A >0$ and big enough $n$, we have  \begin{equation}  \label{tab}  \wt(2\rho-2\rho_{\mathcal{J}})  \geq  A (n \log n \cdot \noJ +  (n-z) n \log(n) ) + \alpha n^2. \end{equation}
for an absolute constant $\alpha$ (we can take $\alpha=1/64$) and for all  $\mathcal{J}$ with $\noJ > 1$. 
\end{Lemma}

Note that the clumsy shape of the right-hand side is just chosen to match with what we will need later.

\proof 
We saw that $\wt(2\rho-2\rho_{\mathcal{J}})$ can be written thus:
$$\frac{1}{2} \sum_{k}  \sum_{a \in J_k} \sum_{b \notin J_k}  |b-a|$$
where the $1/2$ comes from the fact that for each such pair $(a,b)$, either $a>b$ or $b>a$. 

Fixing $k$ for a moment, write $s = n-|J_k|$. 
Note that, for $a\in J_k$, we have 
$$\sum_{b \notin J_k} |a-b| \geq s^2/4$$
since given a set $S \subset \mathbb{Z}$ of size $s$ and not containing $0$, 
we  have $$\sum_S |t| \geq  \begin{cases} (s+1)^2/4-1/4 \\ (s+1)^2/4 \end{cases}$$ according
to whether $s$ is even or odd, with equality attained  e.g. in the cases  $$\begin{cases} S = \{-1, -2, \dots, -s/2\} \cup \{1, 2, \dots, s/2\} \\  \{-1, -2, \dots,- (s-1)/2\} \cup \{1, 2, \dots, -(s+1)/2\}\end{cases} $$

Therefore,
$$\wt(2\rho-2\rho_{\mathcal{J}}) \geq \frac{1}{8}  \underbrace{ \sum_{k} |J_k| (n-|J_k|)^2}_{:=W}.$$

  First of all, all but one part has size $\leq n/2$ and for such parts $ (n-|J_k|)^2 \geq   (n^2/4)$. This immediately leads to the bound
\begin{equation} \label{boundA} \wt(2\rho-2\rho_{\mathcal{J}}) \geq \frac{1}{32} n^2 (\noJ -1).\end{equation}  



  Next, pick any part $J_k$, with $|J_k|=t$.  For $k' \neq k$ we have $n-|J_{k'}| \geq t$, 
and summing over   all such $k'$ gives  a contribution of $ \geq  (n-t)t^2$ to $W$. 
The part  $J_k$ itself contributes $t (n-t)^2$ to $W$; and therefore
$$W \geqslant  (n-t)t^2+t(n-t)^2 = t (n-t) n.$$

If $z$ is the size of the largest part (or more precisely, the size of a fixed part with maximum size), every other part satisfies
 $n-|J_k| \geq n-z$.
 Thus each other part contributes at least $(n-z)^2 |J_k|$, and summing over  the other parts gives  a contribution of at least
$(n-z)^3$ to $W$. On the other hand, the contribution of the  fixed part with size $z$ itself
is $z(n-z)^2$. In total, we get $W \geq n (n-z)^2$.

Averaging the last two bounds gives:
$W \geq \frac{1}{2}n^2 (n-z)$, that is
\begin{equation} \label{boundB} \wt(2\rho-2\rho_{\mathcal{J}}) \geq \frac{1}{16} n^2  (n-z),\end{equation}  
where $z$ continues to be the size of the largest part. 

Adding our \eqref{boundA} and \eqref{boundB}, we get (for $\noJ \geq 2$ and so $n-z \geq 1$): 
 \begin{eqnarray}4 \cdot \wt  &\geq&  \frac{1}{32} n^2 (\noJ-1) + \frac{1}{32} n^2 (n-z) + \frac{1}{32}n^2 (n-z)\\
 &\geq& \frac{1}{64} n^2 \cdot \noJ + \frac{1}{32} n^2 (n-z) + \frac{1}{32}n^2
 \end{eqnarray}
which certainly implies our desired bound.

%
%
%
    
 \qed
  
%

%
%
%

%
 \section{The combinatorics of a block intertwining}
 Here we examine the contribution of a given pair of blocks of the Levi to the constant term of a degenerate Eisenstein series.
 
More specifically, let $M$ be a Levi subgroup. Following our previous notation,    if we write
 $$\mathcal{B} = [N_{t}+1, N_{t+1}], \ \ \mathcal{C} =  [N_r+1, N_{r+1}].$$ 
 for the  $t+1$st and $r+1$st block of the Levi,
 we examine the contribution of $i \in \mathcal{B},  j \in \mathcal{C}$ 
 to the constant term \eqref{Cterm}, and simplify the resulting expression. Roughly speaking we show that the $\zeta$-factors that occur
 in the constant term cancel in a more or less favorable way, so that the larger values tend to be on the bottom.

Let $\nu \in \mathfrak{a}_{M,0}^*$, as in \S \ref{CharM}. 
For short, write
\begin{equation} \label{shorthand} \nu_B = \nu_{t-1}, \ \ \nu_C = \nu_{r-1}.\end{equation}
$$ n_B = n_t, \ \ \ n_C = n_r,$$
$$ \kappa_B = t+1,\ \ \kappa_C = r+1.$$
so these are the $\nu$-values for $\mathcal{B}$ and $\mathcal{C}$, and the size of the blocks $\mathcal{B}$ and $\mathcal{C}$, 
and finally the sequential position of $\mathcal{B}$ and $\mathcal{C}$ respectively. (The $\kappa$-notation will only be used later; we include it here just for reference.) 

  We assume $t< r$, i.e. $\mathcal{B}$ precedes $\mathcal{C}$;  because of our conventions  (see \S \ref{levis}) we have $n_B \leq n_C$.

 \begin{Theorem}  \label{CombinatorialTheorem}Notations as above, let 
 $\upsilon$ be as in \eqref{Mchar}, so that
  {\small  \begin{equation} \label{Mchar}  \upsilon = \left( \frac{n_1-1}{2}+ \nu_1, \frac{n_1-3}{2} +\nu_1, \dots, -\frac{n_1-1}{2} + \nu_1 , \frac{n_2-1}{2} +\nu_2, \dots, - \frac{n_2-1}{2} + \nu_2, \dots \right). \end{equation} }

  Fix $\sigma \in S[M]$, so that $\sigma: \mathcal{B} \coprod \mathcal{C} \rightarrow \{1, \dots, n\}$
 is monotone decreasing on both $\mathcal{B}$ and $\mathcal{C}$ separately. 
 The product  \begin{equation} \label{localproduct}  \prod_{i \in \mathcal{B}, j \in \mathcal{C}, \sigma(i) > \sigma(j)}   \frac{ \xi(\upsilon_{ij}) }{\xi(\upsilon_{ij}+1)}\end{equation} 
(note that for $i \in \mathcal{B}, j \in \mathcal{C}$ we always have $i<j$)  can be rewritten as 
\begin{equation} \label{prod2} \prod_{b \in \mathcal{B}}  \frac{   \xi(\varepsilon_b z + m_b) }{\xi(z+m_b +j_b)}\end{equation} 
where $z = \nu_B-\nu_C, \varepsilon \in \pm 1$ and the $m_b, j_b$ are half-integers.
Moreover, $1/2 \leq m_b \leq n, j_b \geq 0$
and   the  $m_b+j_b$ are pairwise distinct,
and we can suppose   $\varepsilon_b=1$ whenever $m_b=1/2$. 
\end{Theorem}

The (purely combinatorial) proof comprises the rest of this section. 
It will 
be more convenient to index each block by the real parts of the character, i.e. we identify 
 \begin{equation} \label{Aid}\iota_B : \mathcal{B}  \stackrel{\sim}{\longrightarrow} B:= \left\{\frac{n_B-1}{2} , \frac{n_B-3}{2}, \dots, - \frac{n_B-1}{2} \right\}\end{equation}
 via $\iota_B: i \mapsto \frac{n_B-1}{2}  - (N_{t+1}-i)$, and similarly  
\begin{equation} \label{Bid} \iota_C: \mathcal{C} \stackrel{\sim}{\longrightarrow}  C:= \left\{\frac{n_C-1}{2} , \frac{n_C-3}{2}, \dots, - \frac{n_C-1}{2} \right\}\end{equation} 
Observe that the bijections \eqref{Aid} and \eqref{Bid} are order-reversing. 

Finally, set
$$C^* = (C+1/2) \cup (C-1/2) = \left\{ \frac{n_C}{2}  , \dots,  -\frac{n_C}{2}  \right\}$$
(so that $C^*$ and $C$ interlace one another).

We may regard $\sigma$ as a map
$$ B \coprod C \hookrightarrow \{1, \dots, n\}$$
(i.e., the map given by $\sigma \iota_B^{-1}$ on $B$ and $\sigma \iota_C^{-1}$ on $C$). 
Note that $\iota_B, \iota_C$ reverse ordering; in particular, when considered as above, $\sigma$ is monotone {\em increasing}
on $B$ and $C$ individually. 
   With this convention, we can rewrite \eqref{localproduct} as
   \begin{equation} \label{localproduct2}  \prod_{b \in B, c \in C, \sigma(b) > \sigma(c)}   \frac{ \xi(z+b-c) }{\xi(z+b-c+1)}\end{equation}

  Each $b \in B$ specifies a ``cut'' of $C$: it  separates it into the sets given by
  $\sigma(c) > \sigma(b)$ and $\sigma(c) < \sigma(b)$. 
This coincides with the ``cut'' specified by a unique element of $C^*$; said differently, there is a  unique function   
  $f : B \rightarrow C^*$
with the property that for $b \in B, b \in C$ we have
$$\sigma(b) > \sigma(c) \iff f(b) > c.$$

Note that the resulting function $f$ is necessarily non-decreasing because $\sigma$ was increasing on $B$. 
Write $[f(b)] = f(b)-1/2$.  In other words, $[f(b)]$ is the largest element of $C$ that is less than $f(b)$,
if such an element exists;  
if there exists no such element,  which happens exactly when $f(b) = -n_C/2$, we  have $[f(b)] = \frac{-n_C-1}{2}$. 
Therefore, 
$[f(b)]$ is valued in $C \cup \{ -\frac{n_C+1}{2} \}$.

Now fix $b$ and compute $ \prod_{b \in C, \sigma(b) > \sigma(c)}   \frac{ \xi(b-c+z) }{\xi(b-c+z+1)}$:
$$   \frac{ \xi(b-[f(b)]+z) }{\xi(b-[f(b)]+z+1)} \cdot   \frac{ \xi(b-[f(b)]+z+1) }{\xi(b-[f(b)]+z+2)} \dots   \frac{ \xi(b- \left( -\frac{n_C-1}{2} \right)+z) }{\xi(b-\left(- \frac{n_C-1}{2} \right)  + 1+ z )}$$
$$ = 
 \frac{ \xi(b-[f(b)]+z) }{\xi(b +\frac{n_C+1}{2}+z))}$$ 
 (note that this remains valid even in the case when  $f(b) =-n_C/2$ and so $[f(b)] =  \frac{-n_C-1}{2}$:  the product on the left-hand side is empty.)
Introducing a product over $b \in B$, we get 
$$ \prod_{b \in B, c \in C, \sigma(b) > \sigma(c)}   \frac{ \xi(b-c+z) }{\xi(b-c+z+1)} =  \frac{ \prod_{b \in B}  \xi(b-[f(b)]+z) }{\xi(\frac{n_B+n_C}{2}+z) \dots \xi(\frac{n_C-n_B}{2} +1 + z )}$$
 There are the same number $n_B$ of $\xi$-factors in numerator and denominator.  We rewrite it, using the functional equation $\xi(s) = \xi(1-s)$, as 
\begin{equation} \label{rpl2}  \frac{ \prod_{b \in B}  \xi( \frac{1}{2} + |b-[f(b)]-1/2| \pm z) }{\xi(\frac{n_B+n_C}{2}+z) \dots \xi(\frac{n_C-n_B}{2} +1 + z )}=
\frac{ \prod_{b \in B}  \xi( \frac{1}{2} + |b-f(b)| \pm z) }{\xi(\frac{n_B+n_C}{2}+z) \dots \xi(\frac{n_C-n_B}{2} +1 + z )}
\end{equation} 
 where, writing $\varepsilon_b$ for the sign in front of $z$
in the $b$th term of the product (just as  in the theorem statement), 
we have $\varepsilon_b=-1$ precisely when
  $b < f(b)$.

%
%

Recall $n_C\geq  n_B$ because $B$ precedes $C$.     

The following Lemma implies the first statement of the theorem:

\begin{Lemma}  \label{CoreCombinatoriaLemma}
Suppose that $n_C \geq n_B$, and that $$f^*:  B= \{-(n_B-1)/2, \dots, (n_B-1)/2\} \rightarrow \{-n_C/2, \dots, n_C/2\} = C^*$$ is a non-decreasing function.   Then for $q \geq 0$ integer,  the equation \begin{equation} \label{mai} |b-f^*(b)| \geq \frac{n_B+n_C-1}{2} - q\end{equation}  has at most $q+1$ solutions. 
Moreover, equality holds for all $q$ only when either $f^* \equiv -n_C/2$ or $f^* \equiv n_C/2$.
\footnote{These correspond to the cases where every element of $\sigma(B)$ is either less than $\sigma(C)$, or vice versa. 
In the case $f^* \equiv -n_C/2$, the $\xi$-ratio is identically $1$.}

More explicitly, 
 $$ \# \{b : |b-f^*(b)| \geq \frac{n_B+n_C-1}{2} \} \leq 1$$
  $$ \# \{b : |b-f^*(b)| \geq \frac{n_B+n_C-1}{2} -1 \} \leq 2$$
and so on. 

In particular, if we order the the quantities $\frac{1}{2} + |b-f^*(b)|$, for $b \in B$, with multiplicity and in nonincreasing order as
$r_0 \geq r_1 \geq \dots \geq r_{|B|-1}$,
then we have $r_t \leq \frac{n_C+n_B}{2}-t$. 

 \end{Lemma}
 
\proof    
We will prove in the next Lemma that the quantity \eqref{mai} is maximized (for any fixed $q$) by a constant function,
i.e. $f^*$ sending all of $B$ to $c \in C^*$.  Assuming this, we verify \eqref{mai}. By symmetry, we may assume $c \geq 0$. We consider the effect of increasing $c$ by one. Denote by $s_c(t)$ the number of solutions $b \in B$ to $|b-c| = t$. Then $s_c(t) = s_{c+1}(t)$ except when $t = |\frac{n_B-1}{2} - c|$ or $t = |-\frac{n_B-1}{2} - 1 - c|$; in the former $s_{c+1}(t) = s_c(t) - 1$, and in the latter $s_{c+1}(t) = 1$ and $s_{c}(t) = 0$. This shows that, as $c$ increments, the number of solutions to $|b-f^*(b)| \geq r$ does not decrease. This completes the proof of \eqref{mai}.

The remaining assertions follows easily:
 If $r_t > \frac{n_C+n_B}{2}-t$, it follows that 
$\frac{1}{2} + |b-f^*(b)|  \geq \frac{n_C+n_B}{2}-t +1$ 
has at least $t+1$ solutions, i.e. $|b-f^*(b)| \geq \frac{n_C+n_B-1}{2} -(t-1)$
has $\geq t+1$ solutions, contradicting the first assertion. 
\qed

 \begin{Lemma}
 Suppose that $n_C \geq n_B$, and that $$f^*:  B= \{-(n_B-1)/2, \dots, (n_B-1)/2\} \rightarrow \{-n_C/2, \dots, n_C/2\} = C^*$$ is a non-decreasing function. 
Then there exists a constant function $g^*: B \rightarrow C^*$ such that, for any $r$, the number of solutions to $|b-g^*(b)| \geq r$ is no less than the number of solutions to $|b-f^*(b)| \geq r$. 
 \end{Lemma} 
\proof 
Fix $r$. The strategy is to modify $f^*$ one step at a time so that the number of solutions to $|b-f^*(b)| \geq r$ does not decrease after each step \emph{and} $f^*$ becomes constant at the end of the process. 
 
 A {\em modification} of $f^*$ will be a function $g^*$ with the property that $|x-g^*(x)| \geq |x-f^*(x)|$ for all $x \in B$. 
In particular, the number of solutions to $|b-g^*(b)| \geq r$ is at least the number of solutions to $|b-f^*(b)| \geq r$. 

Enumerate elements of $B$ by $b_1, b_2, \ldots,$ in increasing order.  
Suppose that   \begin{equation} \label{f*cons} f^*(b_1) = f^*(b_2) =  \dots = f^*(b_{k-1}) =  f^*(b_k)  <  f^*(b_{k+1}),\end{equation}
where we allow $k=0$ to mean there is no constraint at all.   
%
%
We will show that  either
 \begin{equation} 
\label{caselib} g_1^*(b) = \begin{cases} f^*(b), b \neq b_{k+1} \\ f^*(b_k), b=b_{k+1} \end{cases} \mbox{ or } g_2^*(b) = \begin{cases} f^*(b_{k+1}), b \leq b_{k} \\ f^*(b), b \geq b_{k+1}. \end{cases}.\end{equation} 
is a modification of $f$. Note that $g_j^*(b_1) = \cdots = g_j^*(b_{k+1})$ for $j=1,2$. 
  Iteratively applying this claim shows that we may suppose that $f^*$ is constant, as claimed.

We analyze the following three cases in turn: \begin{enumerate}[(i)]
\item $|b_k - f^*(b_k)| > |b_{k+1} - f^*(b_{k+1})|$.
\item $|b_k - f^*(b_k)| < |b_{k+1} - f^*(b_{k+1})|$.
\item $|b_k - f^*(b_k)| = |b_{k+1} - f^*(b_{k+1})|$.
\end{enumerate} 
\begin{itemize}
\item[(i)]  If $f^*(b_k) > b_k$ then
$$f^*(b_k) - b_k > f^*(b_{k+1}) - b_{k+1} \implies f^*(b_k)+1 > f^*(b_{k+1})$$ which contradicts \eqref{f*cons}.  
Since  $f^*(b_k) \neq  b_k$ by assumption \eqref{f*cons}, we must have  $f^*(b_k) < b_k$; then
 $$ |b_{k+1} - f^*(b_{k+1})| < b_k - f^*(b_k) < b_{k+1} - f^*(b_k),$$
so we can take $g_1^*$ as the desired modification.

 \item[(ii)] If $f^*(b_{k+1}) < b_{k+1}$, then also $f^*(b_k)  < b_k$  by \eqref{f*cons}, and  so
  $$b_k - f^*(b_k) < b_{k+1} - f^*(b_{k+1})   \implies f^*(b_k)+1 > f^*(b_{k+1}),$$ 
 a contradiction as before.
 Therefore  $f^*(b_{k+1}) > b_{k+1}$.  We take $g^*=g_2^*$.
To verify this works, we must check that
for all $b \leq b_k$ we have 
\begin{equation}  \label{oink2} |f^*(b) - b| \leq |f^*(b_{k+1}) - b|.\end{equation}

If $f^*(b_k) \geq b_k$ \eqref{oink2} is true because 
 $$  |f^*(b_{k+1}) - b| = f^*(b_{k+1}) - b \geq f^*(b_k) - b  = f^*(b)  -b \geq 0$$
 
If $f^*(b_k) < b_k$ \eqref{oink2} is true because
$$b - f^*(b) < b_k - f^*(b_k) < f^*(b_{k+1}) - b_{k+1} < f^*(b_{k+1}) - b, $$
$$f^*(b) - b = f^*(b_k) - b < f^*(b_{k+1}) - b.$$
  \end{itemize}

\item[(iii)] We are assuming that $|b_k - f^*(b_k)| = |b_{k+1} - f^*(b_{k+1})|$.
Then either:
\begin{itemize}
\item[(iii-a)] $f^*(b_k) \geq  b_k$. In this case, $g_2^*$ is a modification  just as before -- for $b \leq b_k$, 
 $$   f^*(b_{k+1}) - b \geq f^*(b_k) - b    \geq 0$$  
%

\item[(iii-b)] $f^*(b_k) < b_k$: in this case,  $g_1^*$ is a modification because
 $$   b_{k+1} - f^*(b_k) \geq b_k - f^*(b_k)  = |b_{k+1} -f^*(b_{k+1})|. $$
\end{itemize}

\qed

\section{Bounds}

  \begin{Lemma}  (Bounds for $\zeta$) \label{zetasubconvexity}
 Suppose $|\sigma| \leq 10^{-3}$. Put $q(z) = \min(10|z|, 1)$.   Put $\zeta_q(s)= q(s-1) \zeta(s)$. 
 Then there exists a constant $A$ such that: 
  \begin{itemize}
 \item[(i)] $|\zeta_q(1/2+\sigma+ it)| \leq A (1+|t|)^{0.2}$  
 \item[(ii)] $|\zeta_q(1+\sigma+it)| \leq A (1+|t|)^{0.01}$  
 \item[(iii)] {\em assuming the Riemann hypothesis},  $|\zeta(1+\sigma+it)|^{-1} \leq A (1+|t|)^{0.01}$; 
 \item[(iv)] $|\zeta_q(j+\sigma+it) | \leq A$ and $|\zeta(j+\sigma+it)|^{-1} \leq A$ if $j\geq 3/2$. 
 \end{itemize}
  \end{Lemma}
 \proof For (i), (ii), (iv) see \cite[Chapter 8]{IK}. For  (iii) see Corollary 13.22 of \cite{MV}.
  \qed

 \begin{Lemma} \label{zetabound} (Bounds for $\zeta$ ratios)  
Assume the Riemann hypothesis (this is the only point where we do). 
Let $\zeta_q$ be as in Lemma \ref{zetasubconvexity}. Then  
 uniformly for $|\mathrm{Re}(z)| \leq 0.01/n^4$ we have: 
 $$\left| \frac{\zeta_q(1-z)}{\zeta_q(1+z)}\right| \leq  A \log(5+|t|)(1+|t|)^{|\mathrm{Re}(z)|} \ \ \ t = \mathrm{Im}(z), $$ 
 with an absolute constant $A$ 
\end{Lemma} 

\proof Write $z=\mu+it$. We may as well suppose that $t \geq 10$. 
Split into two cases, according to whether $|\mu| > 1/\log\log(t)$ or not.
We use the bounds in In  Corolllary 13.16 of \cite{MV}:
\begin{itemize}
\item
If $|\mu| < 1/\log\log(t)$, the bounds there show at once that the ratio in question is bounded by an absolute constant times $\log\log(t)$.

 \item
If  $|\mu| > 1/\log\log(t)$,  the bounds therein show that the ratio is bounded by 
$\log(t) \cdot e^{C \log(t)^{2|\mu|}}$ for a suitable absolute constant $C$; without loss $C \geq 2$. Now, 
$$   \log(t) \cdot e^{C \log(t)^{2|\mu|}} \leq A  \log(t) (1+|t|)^{|\mu|},$$
for suitable $A$. \footnote{
  Here we use the fact that, with $u=\log(t)$, we have the bound 
  $$ Cu^{2\mu} \leq  C^4+ \mu u, u \geq C^4$$
To verify this, look at the function $C^4+ \mu u- C u ^{2 \mu}$. Its derivative is $\mu(1-2 C u^{2 \mu-1})$ , so it is increasing for $u \geq C^2$.}
\end{itemize}

 \qed 
 \begin{Lemma} \label{Gammabound} (Bounds for $\Gamma$) 
 
    Write, for short, $\Gamma_{\R}(s) = \pi^{-s/2} \Gamma(s/2)$. 
  We have the bounds, valid uniformly for $\mathrm{Re}(z) \geq 0.49$ and $|p| \leq 0.26$ (with an absolute constant $A$ and $R$): 
    \begin{equation}   \label{gammabound0}  A^{-1} \leq \frac{1}{|z|^{p/2}}  \left| \frac{\Gamma_{\R}(z+p)}{\Gamma_{\R}(z)} \right|  \leq A
   \end{equation}   \begin{equation}   \label{gammabound1}   \left| \frac{\Gamma_{\R}(z+r)}{\Gamma_{\R}(z)}\right| \geq (1+|z|)^{1/4}  \ \ (r \in  \frac{1}{2} \mathbf{Z}, r  > R). \end{equation}
 

 \end{Lemma}
 \proof 
Recall Stirling's formula (see, for example, 6.1.40 of \cite{AS})
 $$\Gamma(z) = \sqrt{2\pi} z^{z-1/2}e^{-z} (1+ \varepsilon)$$
 where $|\varepsilon| \leq  \frac{1}{4 |z|}$; this formula is valid for $\mathrm{arg}(z)$ between $(-\pi/2, \pi/2)$. 
 
 For \eqref{gammabound0} it is enough to prove a similar bound with $\Gamma_{\R}$ replaced by $\Gamma$, i.e. upper and lower bounds
 for $ \frac{1}{|z|^{p}}  \left| \frac{\Gamma(z+p)}{\Gamma(z)} \right|  $      where we now restrict to $\mathrm{Re}(z) \leq 0.49/2$ and $|p| \leq 0.49/2$. 
Stirling says that $\left| \frac{\Gamma(z+p}{\Gamma(z)} \right| =  \frac{z^{p}}{e^{1/4}}   \left(1 + \frac{p}{z} \right)^{z+p-1/2}  \cdot  \left( \frac{1+\varepsilon_1}{1+\varepsilon_2} \right)$.
 Now the bracketed term is absolutely bounded in the specified region.
 So is  both
 $$  \log(1+\frac{p}{z} ) \mbox{ and }  \log  \left(1 + \frac{p}{z} \right)^{z} = z  \log(1+\frac{p}{z}) $$
(in the first case $z \mapsto  p/z$ takes
the set $\mathrm{Re}(z) > 0.49$  into a subset of the disc $|u| \leq 0.9$; in the second case, we use the power series for $\log(1+p/z)$ when $z$ is large.) Our claim \eqref{gammabound0} follows. 
 
 The claim \eqref{gammabound1} follows from \eqref{gammabound0} applied with $p=1/2$ many times, once we notice
 $$A^{-N-1} |z (z+1) \dots (z+N)| \geq  \frac{ N! A^{-N-1} }{3} \geq 1$$
 for an $N$ that depends only on $A$. 
%
%
%
 \qed 
 
\begin{Theorem}
Let $B$ and   $j_b, m_b$ be as in Theorem \ref{CombinatorialTheorem}, that is to say:
 $1/2 \leq m_b \leq n, j_b \geq 0$, the $m_b+j_b$ are pairwise distinct,
and     $\varepsilon_b=1$ whenever $m_b=1/2$. 

Write  $q(z) = \min(10 |z|,1)$. We have a bound
  valid uniformly for $|\mathrm{Re}(z)| \leq 0.001/n^4$, with an absolute constant $C$:
\begin{equation} \label{frobound}    \prod_{B}  \frac{      \xi(\varepsilon_b z + m_b) }{\xi(z+m_b +j_b)} \leq C^{|B|} \prod_{b: m_b = 1}q(z)^{-1}
\cdot  \log(5+|t|) \cdot  \begin{cases}  (1+|t|)^{-0.02}, \mbox{ some $j_b \neq 0$,} \\
(1+|t|)^{0.01/n^2}, \mbox{ all $j_b=0$.} \end{cases} 
\end{equation}  Finally,    when all   $j_b= 0$, the corresponding product has absolute value $1$ on the line $\mathrm{Re}(z) = 0$ (obvious). 
\end{Theorem}
\proof 
Set $\xi_q(s) = \pi^{-s/2} \Gamma(s/2) \zeta(s) q(s-1)$.  Note that, for $|\mathrm{Re}(z)| \leq 0.01$, 
we have $\xi_q(z+m_b) = \xi(z+m_b)$ unless $m_b=1$. 
What we will prove is that (for a suitable absolute constant $C$):
\begin{equation} \label{subproof} \left|   \frac{       \xi_q(\varepsilon_b z + m_b) }{\xi(z+m_b +j_b)}   \right| 
\leq C   \begin{cases} 
 (1+|t|)^{-0.03},    &  \mbox{ unless $j_b=0$},   \\ 
 (1+|t|)^{0.01/n^4},  &  \mbox{   $j_b$} = 0, m_b \neq 1 \\
   (1+|t|)^{0.01/n^4} \log(5+|t|), & \mbox{ $j_b$}= 0, m_b = 1\\
   \end{cases}  \end{equation} 
The result follows: if we take the product of all these factors, we get at least one factor of $(1+|t|)^{-0.03}$ as long as not all $j$s are zero;
   we get at most $|B| \leq n$ factors of $(1+|t|)^{0.01/n^4}$, and finally at most one factor of $\log(5+|t|)$.

Now to prove \eqref{subproof}.  We write the quantity to be bounded as
\begin{equation} \label{oinkoink}  \frac{\zeta_q(\varepsilon_b z+m_b)}{\zeta(z +m_b+j_b)} 
 \cdot \frac{\Gamma_{\R}(\varepsilon_b z+m_b)}{\Gamma_{\R}(z +m_b+j_b)}  \end{equation}
 We subdivide into four cases. The first three cases will be $j_b > 0$ but $j_b$ in different ranges;
 the last case is $j_b=0$.

 \begin{itemize}
 \item  First case of $j_b \neq 0$: 
 $j_b \geq  R+1$, where $R$ is the absolute constant in the Lemma \ref{Gammabound}.

  By Lemma \ref{zetasubconvexity} the  $\zeta$ in the denominator is absolutely bounded below.   
By the same Lemma, the $\zeta_q$ in the numerator is in all cases bounded by a constant multiple of $(1+|t|)^{0.2}$. 
Therefore,   $\left| \frac{\zeta_q(z)}{\zeta(z+m_b+j_b)} \right| \leq C (1+|t|)^{0.2}$. 
 
If $\varepsilon = 1$, we apply the bounds  \eqref{gammabound1} 
directly, with $r=j_b$, to get that the $\Gamma$-ratio is bounded by  $C(1+|t|)^{-0.25}$; we are done.

If $\varepsilon  = -1$, we write $z= \sigma+it$ and conjugate the numerator to get: 
$$ \left| \frac{\Gamma_{\R}(\varepsilon_b z+m_b)}{\Gamma_{\R}(z +m_b+j_b)} \right| =
\left| \frac{\Gamma_{\R}(-\sigma+it+m_b)}{\Gamma_{\R}(\sigma+it +m_b+j_b)} \right| \leq C (1+|t|)^{-0.24}$$
where we applied  \eqref{gammabound0} with $r=2\sigma$ and then \eqref{gammabound1} with $r=j_b$ to conclude. 

\item  Second case of $j_b \neq 0$:  suppose that $0< j_b \leq R$ but $m_b \geq R$, where
$R$ is the absolute constant in Lemma \ref{Gammabound}. Referring again to \eqref{oinkoink}, 
 the $\zeta$-quotient is absolutely bounded above, and iterated application of \eqref{gammabound0} at most $4R$ times (so the implicit constants don't matter) 
 gives the bound $C (1+|t|)^{-0.24}$ for the ratio of $\Gamma$-functions. 
 This proves the desired bound this case.

 \item  Third and final case of $j_b \neq 0$:  $0 <  j_b \leq  R$  and $m_b <  R$.
 
 There are at most $O(1)$ of these factors,
 because each factor of the denominator occurs $O(1)$ times --  recall that, by assumption, the various $m_b+j_b$ are pairwise distinct. 
   We can therefore ignore any implicit constants in this analysis.  \  So fix $m,j \in \frac{1}{2} \mathbf{Z} \cap [1/2, R]$ and  examine again  \eqref{oinkoink}.
 The $\Gamma$-quotient decays as $(1+|t|)^{-0.24}$ by the same logic as in the previous case, and the quotient of $\zeta$-functions grows at most as $(1+|t|)^{0.21}$ by Lemma \ref{zetasubconvexity}. This proves the desired bound in this case.

\item Terms with $j_b= 0$.  Here we can assume that $\varepsilon_b = -1$
(otherwise the term is $1$) and therefore also that $m_b \geq 1$ (part of our assumptions: see last sentence of Theorem statement). 

The term in question is
$$   \frac{\zeta_q(- z+m_b)}{\zeta(z +m_b)} 
 \cdot \frac{\Gamma_{\R}(- z+m_b)}{\Gamma_{\R}(z +m_b)}  $$
 
The $\Gamma$-ratio is  bounded above, by means of \eqref{gammabound0},   by an absolute constant 
multiplied by $(n+|t|)^{0.001/n^4} \leq C(1+|t|)^{0.001/n^4}$ (here we also used the fact that $m_b \leq n$). 
 If $m_b > 1$, the $\zeta$-term is absolutely bounded above and below. 
 If $m_b=1$ (and note that there is at most one term of this form, in the product we are analyzing)
  the $\zeta$-term $\zeta_q(-z+1)/\zeta_q(z+1)$ can be analyzed with Lemma \ref{zetabound} to get a bound of  $\log(5+|t|)  \cdot (1+|t|)^{0.001/n^4}$.

  \end{itemize}
  
  This concludes the proof of the Theorem. 
 
   \section{ Bounding the constant term via a contour integral}
 
We now return
to analyzing the behavior of the full constant term  \eqref{Cterm} for $E_{M,\nu}$. As mentioned after \eqref{Cterm}, some of the individual terms in \eqref{Cterm} can have poles;
however,   $(E_{M,\nu})_N$ and so the whole sum is holomorphic at the points  of interest where $\mathrm{Re}(\nu_i)= 0$. 
To bound
\eqref{Cterm}, then, we deform along a contour where none of the individual terms have poles, and use the Cauchy estimate. 
 
 \begin{Theorem} \label{secondth} Let $f: A \rightarrow \C$ be as in \eqref{fdef}. Let $M$ be a Levi subgroup with $k$ blocks, and let $\nu = (\nu_1, \dots, \nu_k) \in \mathfrak{a}_{M,0}^*$ parameterizes a unitary character of $M$.
 Let other notation 
 be as in \S \ref{sec:degenerate}.   Then  for an absolute constant $C$, 
 \begin{equation} \label{tbb}   |  \langle f a^{-2 \rho},  (E_{M,\nu})_N\rangle_A  |
\leq  C^{nk \log(n)} \sum_{\sigma \in S[M]}    | \langle f, a^{-\rho+ \sigma \upsilon} \rangle_A|  \end{equation}
where $k$ is the number of parts of the partition associated to $M$, $\upsilon$ is the shifted $\nu$-parameter,
as in \eqref{Mchar}, and $S[M]$ is as in \eqref{SMdef}.  The inner products on both sides are taken in $L^2(A)$. 
\end{Theorem}

\proof

Note that the right-hand side of \eqref{tbb} equals, by \eqref{frofro},  $ \sum_{\sigma}\left|  \frac{ T^{\mathrm{Re}(\wt(\rho-\sigma \upsilon))}  C^{n  k \log(n)}}{ \prod |\mu_i(-\rho + \sigma \upsilon)|}  \right| $.

The left hand side of \eqref{tbb} can be rewritten, following  \eqref{constant2}, and Theorem  \ref{CombinatorialTheorem}, as  $2^{1-n}$ times the absolute value of
\begin{equation} \label{tbb2}  \sum_{ \sigma \in S[M]} \prod_{B< C}  \prod_{b \in B}  \frac{ \xi(m^{BC\sigma}_b +\varepsilon_b^{BC\sigma}(\nu_B-\nu_C)) }{\xi(m^{BC\sigma}_b +j^{BC\sigma}_b + \nu_B-\nu_C)}    
\langle f, a^{-\rho+\sigma \upsilon} \rangle \end{equation} 
 where   $m^{BC\sigma}, j^{BC\sigma} \in \frac{1}{2}\mathbf{Z}_{\geq 0}$  and $\varepsilon^{BC\sigma} \in \{\pm 1\}$ are so denoted to recall that they depend on $B, C$ and $\sigma$. 
 Moreover we wrote for short $\nu_B = \nu_r$ if $B$ is the $(r+1)$st block $[N_{r}+1, \dots, N_{r+1}]$, just as in \eqref{shorthand}.
 For each fixed $\sigma, B, C$ the   $m^{BC\sigma}, j^{BC\sigma} , \varepsilon^{BC\sigma} \in \{\pm 1\}$ 
 satisfy the constraints enunciated in Theorem \ref{CombinatorialTheorem}.

We will bound this  left-hand side one $\sigma$ at a time. 

First of all, let us separate into cases according to whether $j^{BC\sigma} = 0$ for all $B, C$ or not.

{\em Case 1:}  $\sigma$ is such that  $j^{BC\sigma} = 0$ for all $B, C$. 
These terms are trivially bounded by $\langle f, a^{-\rho+ \sigma \upsilon} \rangle_A$: 
the term looks like  
$$  \prod_{B,C}  \prod_{b \in B}  \frac{ \xi(m^{BC\sigma}_b +\varepsilon_b^{BC\sigma}(\nu_B-\nu_C)) }{\xi(m^{BC\sigma}_b  + \nu_B-\nu_C)} \langle f, a^{-\rho+\sigma \upsilon} \rangle$$ 
The function $\xi(m-z)/\xi(m+z)$ is holomorphic everywhere along the line $\mathrm{Re}(z) = 0$. Therefore, 
the $\xi$-product/ratio is holomorphic at the given value of $\nu_B$, and even better,  has absolute value $1$.

 {\em Case 2:} there exists at least one $B,C$ for which $j^{BC\sigma} > 0$. 
To bound this, we apply \eqref{frobound} and a contour integration argument.   Fix some small parameters 
$$ \delta =n^{-8}, \  \ \ H=n^{-4}.$$
 
As in \eqref{shorthand}, let $\kappa_B$ be the parameter of the block $B$ (i.e., the block $[1, N_1]$ has parameter $1$ and so on). 

  Let $\mathcal{C}_{\delta,H}$ be the {\em basic contour}: the oriented closed curve in the complex plane 
which consists of the two vertical lines $\pm \delta + i v$, for $v \in [-H, H]$,
together with the two horizontal lines $\pm i H + b, b \in [-\delta, \delta]$.   

Let $\nu_z \in \mathfrak{a}_M^*$ be the parameter:
$$ \nu_z := (\nu_1+z, \nu_2 + 2z, \dots, \nu_k + kz)$$
so that $\nu_z$ takes the value $\nu_B + \kappa_B z$ on the block $B$.  Note that $\nu_z$ will no longer be unitary,
i.e. $\nu_z \notin \mathfrak{a}_{M,0}^*$. 
We define $\upsilon_z$ to be the shifted $\upsilon$-parameter attached to $\nu_z$,
so that the relationship between $\upsilon_z$ and $\nu_z$ is as in \eqref{Mchar}.

 Note that $E_{M,\nu_z}$ is a meromorphic function of $z$,
and we will use Cauchy's formula to compute $E_{M,\nu}$ as a contour integral along $\mathcal{C}_{\delta,H}$.
As $z$ moves along $\mathcal{C}_{\delta,H}$, 
  the term $\nu_B-\nu_C$ moves along the contour $\nu_B-\nu_C + (\kappa_B-\kappa_C) \mathcal{C}_{\delta,H}$.
  Note that $\kappa_B -\kappa_C$ is nonzero for $B \neq C$, so $z$ is really moving. 
  
Thus we must study \eqref{tbb2}  as $z \in \mathcal{C}_{\delta,H}$.   
For $\nu$ of the form $\nu_z$, with $z$ on $\mathcal{C}_{\delta,H}$, we have by \eqref{frobound}:
 \begin{eqnarray*} \left| \prod_{B<C}  \prod_{b \in B}  \frac{ \xi(m^{BC\sigma}_b +\varepsilon_b^{BC\sigma}(\nu_{B,z}-\nu_{C,z})) }{\xi(m^{BC\sigma}_b +j^{BC\sigma}_b + \nu_{B,z}-\nu_{C,z})}     
 \right|  &\leq&    \prod_{B<C}   q(    \underbrace{\left( \nu_B-\nu_C + (\kappa_B-\kappa_C) z\right)}_{\nu_{B,z}-\nu_{C,z}} )^{-r(B, C)} \\
\log(5+|t|)^{k^2}  \times \prod_{B,C} C^{|B|} &  \times &   (1+|t|)^{0.01 \aleph/n^2 - 0. 02 \beth}\end{eqnarray*} 
where $t = \mathrm{Im}(\nu_{B,z}-\nu_{C,z})$, and:
\begin{itemize}
\item  $\aleph$ is the number of $(B,C)$ where all $j^{BC\sigma}_b = 0$  
\item  $\beth$ is the number of $(B,C)$ where  not all $j^{BC \sigma}_b= 0$. 
\item $r(B,C)$ is   the number of $b \in B$ with $m^{BC \sigma}_b = 1$. 

\end{itemize}
By assumption, $\beth \geq 1$
and of course $\aleph \leq n^2$; thus $(1+|t|)^{0.01 \aleph/n^2 - 0. 02 \beth} \leq (1+|t|)^{-0.01}$.   Next
$$\max_t \log(5+|t|)^{k^2} (1+|t|)^{-0.01} \leq C^{\log(k) k^2}.$$
Also, $\prod_{B,C} C^{|B|} \leq C^{nk}$, 
so the above expression is actually bounded by $$ C^{nk \log(n) }  \prod_{B,C}   q(   \left( \nu_B-\nu_C + (\kappa_B-\kappa_C) z\right) )^{-r(B, C)} .$$

So it remains to give an upper bound for \begin{equation} \label{georr} \int_{\mathcal{C}_{\delta,H}} \prod_{B,C}   q(   \left( \nu_B-\nu_C + (\kappa_B-\kappa_C) z\right) )^{-r(B, C)}\end{equation} 
(a priori this quantity could even be infinite for certain $\delta, H$).

Recall that $q(z) = \min(10 |z|,1)$ and in particular
$$q(z)^{-1} \leq \min\left(|\mathrm{Re}(z)|^{-1}, |\mathrm{Im}(z)|^{-1}\right).$$

Write
$t_B =\mathrm{Im}(\nu_B), t_C = \mathrm{Im}(\nu_C).$
Then for any $z \in \mathcal{C}_{\delta,H}$, by the inequality just above, we have  
\begin{eqnarray}\label{ribbitribbit}  \ \ q(\nu_B-\nu_C+(\kappa_B-\kappa_C) z)^{-1}  \leq \delta^{-1} &+&   \frac{1}{|(t_C-t_B)  -  (\kappa_C-\kappa_B) H|} \\ &+&  \frac{1}{|(t_C-t_B)  +  (\kappa_C-\kappa_B) H|. \nonumber }\end{eqnarray}

\begin{Lemma} Write $\log^+(x) = \max(\log(x), 0)$ for $x > 0$, and $\log^+(x) =0$ for $x \leq 0$. 
 Then:
 \begin{itemize} \item For $x,y> 0$, we have  $\log^+(x+y) \leq \log^+(x) + \log^+(y) +1$.
 \item  On any subinterval $I \subset \R$ of length $s$, the average
 value of $\log^+(|x|^{-1})$  is  bounded  above by  
$$\begin{cases} 1- \log(s/2), s < 2 \\ 2/s, s \geq 2 \end{cases}.$$ 
\end{itemize} 
\end{Lemma}

\proof 
 Note that $\log^+(x+1) \leq \log^+(x)+1$ for $x \geq 0$: obvious for $x \leq 1$, and otherwise
 it follows from $(x+1) \leq e x$. Thus the result follows for $y \leq 1$ or (symmetrically) $x\leq 1$.  Otherwise,
we must check $\log(x+y) \leq \log(x) + \log(y)+1$, i.e. $x+y \leq exy$; 
without loss $x \geq y$, and then $exy \geq ex \geq 2x \geq x+y$.   

Given any interval $I = [a,b]$ with $0< a < b$, visibly the average value of $\log^+(x^{-1})$
on $I$ is less than its average value on the shifted interval $I-a = [0,b-a]$.   
Next, given any interval $I$ containing zero, the average value of $\log^+(|x|^{-1})$ 
is less than the corresponding average on the ``symmetrized'' interval $[-\ell/2, \ell/2]$,
where $\ell$ is the length of $I$.  

We are reduced to doing the computation on the interval $(0, s/2)$; 
in that case we get (with $b=s/2$, and assuming $b \leq 1$):
$$ \frac{1}{b} \int_0^b -\log(x) = \frac{1}{b}  (-x \log(x) + x)|^b_0 = -\log(b)+1$$

Therefore, the average is  $-\log(s/2) + 1$ if $s/2 < 1$.
If $s/2 \geq 1$ we similarly get an upper bound of $2/s$.
 \qed

Using the obvious bound $\log(\ldots) \leq \log^+(\ldots)$, and taking into account that $r(B, C) \leq |B|$ and so $\sum_{B,C} r(B, C) \leq kn$, we get that for all $z \in \mathcal{C}_{\delta,H}$ 
\begin{eqnarray}
\log \prod_{B,C}   q(\nu_B-\nu_C + (\kappa_B-\kappa_C) z)^{r(B, C)} \leq    C kn \log(n)  \\  \nonumber  +   \sum_{B,C} r(B, C) \log^+(|(t_C-t_B)  -  (\kappa_C-\kappa_B) H| ^{-1}) + \left( \mbox{ similar term} \right).\end{eqnarray}

Now  we  choose $H$ suitably.    Average over $0 < H < 1/n$;  the average of each $\log^+$ term is, by the prior lemma, 
at most $O(\log n)$. 
Therefore, there exists $H \in (0, 1/n)$
for which the right-hand side is bounded above by $C nk \log(n)$. 

In summary, we can choose $H \in (0, 1/n)$ with the property that 
the integrand of \eqref{georr} is bounded  above by $\exp(C n k \log(n))$. Cauchy's integration formula now concludes the proof (of the Theorem \ref{secondth}).  \qed

 \section{Conclusion}
 
 We now complete the paper by giving the proof of (ii) of the main Theorem 
\ref{MainTheorem}.  At this point this is a matter of  putting our prior bounds together, plus some elementary estimates for integrals on Euclidean spaces. 

We continue to use $\langle -, - \rangle_A$ to denote inner products in $L^2(A)$. 
In our notation, (ii) of the main Theorem amounts to the assertion that
 \begin{equation} \label{D1} \|E_f - \bar{E_f}\|^2 < \exp(-a n^2) \|\bar{E}_f\|^2\end{equation} for some $a> 0$ and big enough $n$.
  
  Terms of size $n^n$ are negligible from the point of view of proving \eqref{D1}.
So we shall use the notation 
  $$B \sim_n B'$$
 if there exists constants $b_1, a_2$ so that
 $$ \frac{B}{B'}, \frac{B'}{B} \leq b_1 \exp(a_2 n \log(n))$$
 Often we will supress the subscript $\sim_n$ and just write $\sim$.

%
%
%

 
 \subsection{}

 Begin with the spectral decomposition \eqref{SpectralDecomposition}, or what we obtain from it by taking inner products with $E_f$:
 \begin{eqnarray*} \|E_f - \bar{E}_f\|^2 &=&  \sum_{M \neq G} V_{M}^{-1} \int_{\nu \in \mathfrak{a}_{M,0}^*/W_M}  \left| \langle E_f, E_{M,\nu} \rangle \right|^2  
 \\ &\stackrel{\eqref{unfold2}}{=}& 4^{1-n} \sum_{M \neq G} V_M^{-1}  \int_{\nu} \left| \langle f a^{-2 \rho}, (E_{M,\nu})_N \rangle_A  \right|^2 d\nu  
\\ &\stackrel{\eqref{tbb}}{\leqslant} &   A^{n k \log(n)} \sum_{M}  V_{M}^{-1} \int_{\nu} \left( \sum_{\sigma \in S_M}  \left| \langle f,  a^{-\rho + \sigma \upsilon} \rangle_A \right| \right)^2 
\end{eqnarray*} 
Note that the total number of $M$ and also the size of $|S[M]|$ are both $\sim_n 1$.
Also recall  from \eqref{unfold} that 
$$ \bar{E}_f = V_G^{-1} 2^{1-n}  \langle f, a^{-2 \rho} \rangle_A.$$ 
Therefore, \eqref{D1} will follow if we check that, for any 
constant $A$ and sufficiently large $n$,  there is a constant $A'$ such that 
\begin{equation} \label{Desideratum} A^{nk \log(n)} \frac{ V_M^{-2}  \int_{\nu} |\langle f , a^{-\rho + \sigma \upsilon} \rangle_A |^2 d\nu}{   V_G^{-2}  |\langle f , a^{-2\rho} \rangle_A |^2} \leq A'  \exp(-\delta n^2) \end{equation} 
where $k$ is the number of parts in  the partition associated to $M$, or equivalently in the division $\mathcal{J}$ (see \S \ref{sec:reindexing}) associated to the pair $(M, \sigma)$. 
 
From \eqref{frofro}  we have $|\langle f, a^{-2 \rho} \rangle| \sim_n T^{\wt(2 \rho)}$. Let $s \leq n-1$ be the size of the largest part of $\mathcal{J}$.
By an elementary estimate from the definitions \eqref{Mvoldef} and \eqref{Qdef},   we see 
that $\frac{V_G}{V_M}$ is bounded (up to factors of size $e^{C n \log(n)}$) by $ \xi(s+1) \dots \xi(n)  \leq \mathrm{const} \cdot \xi(n)^{n-s+1}$, 
which leads to the bound:
$$ \frac{V_G}{V_M} \leq  \exp(C  (n-s) n \log(n))$$

 We will prove that for $\sigma \in S[M]$ so that $(M, \sigma)$ is parameterized by $\mathcal{J}$ we have 
\begin{equation} \label{Mellinbound}  \int_{\nu} |\langle f , a^{-\rho + \sigma \upsilon} \rangle_A |^2 \ll_n  C^{n \log n} T^{2 \wt(\rho+\rho_{\mathcal{J}})}\end{equation}
Here  $\mathcal{J}$ is as in \S \ref{rhoJdef}. 
Once that's done,  our conclusion \eqref{Desideratum} follows from Lemma \ref{weightbound}, 
since $\langle f, a^{-2 \rho} \rangle_A^2 \sim_n T^{4 \wt(\rho)}$.

\proof (of the bound \eqref{Mellinbound}):

Recall first that \eqref{realpart} shows that
$$\mathrm{Re}(\sigma \upsilon) = - \rho_{\mathcal{J}}$$
Now, \eqref{frofro} shows that 
$$   |\langle f , a^{-\rho + \sigma \upsilon} \rangle_A |^2 = T^{2\wt(\rho)+2\wt(\rho_{\mathcal{J}}) } \cdot \prod_{i} |\mu_i(-\rho+\sigma \upsilon)|^{-2}$$

The real part of $\mu_i(-\rho+\sigma\upsilon)$, i.e. of $\mu_i(-\rho-\rho_{\mathcal{J}})$,  satisfies
$$| \mathrm{Re}\  \mu_i(\rho+\rho_{\mathcal{J}}) |\geq 1,$$
since both $\rho$ and $\rho_{\mathcal{J}}$ is a sum of positive roots, and each positive simple root occurs at least once. 

For $1 \leq i \leq n-1$, we write   $u_i = \mathrm{Im} \ (\mu_i(\sigma \upsilon))$. 
Therefore, $|\mu_i(-\rho+\sigma \upsilon)|^{-2} \leq \frac{1}{1+u_i^2}$, and it will be enough 
  to bound $\int_{\nu \in \mathfrak{a}_{M,0}^*} \prod_{i} \frac{1}{1+u_i^2}$, 
 
Write  $t_i  = \mathrm{Im}(\nu_i)$;
recall that the subset of unitary characters in $\mathfrak{a}_M^*$  is identified,
 when considered as characters on $Z_M$ (cf. \eqref{Zmchar}) 
 as the subspace  $\subspace \subset \mathbf{R}^n$ where 
$$\subspace = \{\mathbf{t} \in \R^n: t_1=\dots = t_{N_1}, \ \ t_{N_1+1} = \dots = t_{N_2}, \dots  \}$$
and moreover $\sum t_i = 0$.    

Since $\nu$ and $\upsilon$ have the same imaginary part, coordinate by coordinate, we see
\begin{equation} \label{usimple} u_{i} = \mathrm{Im}(\mu_{i} (\sigma \upsilon)) = \mathrm{Im} (\mu_{i}  (\sigma  \nu)) \stackrel{\eqref{reverse_formula}}{=} 
\sum_{j=1}^{i} \mathrm{Im} (\sigma \nu)_j = \sum_{j=1}^{i} t_{\sigma^{-1} j}
  \end{equation} 
  
we wish to prove   $$ \int_{\subspace} \frac{1}{1+t_{\sigma^{-1}(1)}^2} \frac{1}{1+t_{\sigma^{-1}(1)}^2+ t_{\sigma^{-1}(2)}^2} \dots 
\ll_n  C^{n \log n} $$
The measure on $Q$ has been described in \eqref{numeasure};
it differs  by exponentially bounded factors (negligible for our purpose) 
from the measure $\left| \wedge_i d(t_{N_i}-t_{N_i+1}) \right| $ on $\subspace$.
In what follows, we equip $Q$ with this latter measure.

%
%

 
 \subsection{Some Fourier analysis}
 
The proof of \eqref{Mellinbound} now reduces to an elementary estimate on Euclidean spaces (we take a long time to do it just because we get very nervous about measures). 

 Let $V$ be the full space $\{(t_1, \dots, t_n) \in \mathbf{R}^n: \sum t_i = 0\}$, equipped with the measure $$\nu_b= | \wedge_i (dt_i- dt_{i+1}) |.$$ 
 Thus $\dim V = n-1$ and $\subspace \subset V$ has dimension $k-1$;
 the measure on $\subspace$ is as just described.

   For a nice function $f$ on $V$, we have by the usual Fourier inversion formula  
\begin{equation} \label{naive}  \int_\subspace f =   \int_{(\kappa_1, \kappa_2, \dots, \kappa_{n-k}) }  d\kappa \int_{V} f(t)  e^{2 \pi i \left( \kappa_1(t_2-t_1) + \kappa_2(t_3-t_2) + \dots\right) } d\nu_b\end{equation} 
  {\em where we omit the terms corresponding to $t_{N_k+1}-t_{N_k}$}, and the $\kappa$-measure is usual Lebesgue measure. Now the sum in the exponential is given by 
\begin{eqnarray} \nonumber
\kappa_1(t_2-t_1) + \kappa_2 (t_3-t_2) + \dots   + \kappa_{N_1-1}  (t_{N_1}-t_{N_1-1})  + \kappa_{N_1} (t_{N_1+2}-t_{N_1+1}) + \dots \\  = -\kappa_1 t_1 + (\kappa_1-\kappa_2) t_2 +   \cdots
  +  (\kappa_{N_1-2}-\kappa_{N_1-1}) t_{N_1-1} + \kappa_{N_1-1} t_{N_1} + \dots \end{eqnarray} 
 Changing coordinates
 in the $\kappa$ variable,
we can re-write \eqref{naive} as: 
\begin{equation} \label{Finvesion}  \int_\subspace f =   \int_{\subspace^{\perp}} \widehat{f}(k), \ \ \ \widehat{f} := \int_V  d\nu_b(t) \ f(t) e^{2 \pi i \sum k_i t_i} \end{equation} 
    where  $Q^{\perp}$ is the $n-k$ dimensional space
    \begin{equation} \label{Wperdef} \subspace^{\perp} = \{ (k_1, \dots, k_n):  \sum_{N_j+1}^{N_{j+1}} k_i = 0 \}\end{equation} 
    and the measure on $Q^{\perp}$ is that obtained by taking the product of $$dk_1 \wedge \dots \wedge dk_{N_1-1}, \ \    dk_{N_1+1} \wedge \dots \wedge dk_{N_2-1},$$ i.e. omitting one coordinate from each block.

Write as above $u_i =  \sum_{j=1}^i t_{\sigma^{-1}j}$
and $f= \prod_{i=1}^n  \frac{1}{1+u_i^2}$.
  We evaluate $\int_\subspace f$ by means of Fourier inversion \eqref{Finvesion}. 
We have 
 $$  \hat{f}(k) =   \int  d\nu_b \ e^{2 \pi i  \sum k_i t_i} \prod \frac{1}{1+u_i^2}.$$
 For short, let us write $k_i^* = k_{\sigma^{-1}(i)}$ and
 $t_i^* = t_{\sigma^{-1}(i)}$. 
Note that 
 $$\sum k_i t_i =  \sum k^*_i t^*_i
 $$
 {\small $$ = k_n^* (t_1^* + \dots + t_n^*) + (k_{n-1}^*-k_n^*)(t_1^* + \dots + t_{n-1}^*) + (k_{n-2}^*-k_{n-1}^*) (t_1^* + \dots + t_{n-2}^*) + \dots +   (k_1^* - k_2^*) t_1^*$$}
 
 Therefore, if we rewrite the integral for $\hat{f}(k)$ in coordinates $u_{n-1} = t_1^* + \dots t_{n-1}^*, u_{n-2} = t_1^* + \dots + t_{n-2}^*, \dots, t_1^*$, 
 and note that
 $$|du_1 \wedge \dots \wedge du_{n-1} | =| dt_1^* \wedge \dots \wedge dt_{n-1}^*| \stackrel{\S \ref{3meas}}{=}  \frac{1}{n} d\nu_b(t)$$
 we get   \footnote{Recall that $\int e^{2 \pi i k x}/(1+x^2) $ can be evaluated by contour integration to be
 $\pi e^{-2 \pi |k|}$. } 
  $$ \hat{f}(k) =  n \pi^{n-1}  e^{- 2 \pi \|k^*\|},$$
  
 where we write
 $$ \|k^*\| = |k^*_1-k^*_2| + |k^*_2-k^*_3| + |k^*_3-k^*_4| + \dots  + |k_{n-1}^* - k_n^*| + |k_n^*|.$$ 

%
Therefore, applying \eqref{Finvesion}, we get 
$$  \int_\subspace f   \leq   n \pi^{n-1}  \int_{\subspace^{\perp}}  e^{-2 \pi \|k^*\|}  $$
and it remains to check that $\int_{\subspace^{\perp}} e^{-2 \pi \|k^*\|}$ 
is exponentially bounded. 
Let $\Xi$ be the measure of the set $\|k^*\| \leq 1$ inside $\subspace^{\perp}$. Then, by a homogeneity argument, the last integral equals
$  \Xi  \cdot \int e^{-2 \pi z} d(z^{n-k}) =     (n-k)!  (2\pi)^{-(n-k)} \cdot \Xi$. Finally we have in fact $\Xi \leq 1$: 
we easily see that $\|k^*\| \geq \sup_i |k|$. 
Thus it is enough to bound the volume of the set $\sup_i |k_i| \leq 1$ inside $\subspace$.  
But this is a product of the volume  of various sets of the type
inside $$ \{  (y_1, \dots, y_d) \in \mathbf{R}^d: \sup_i |y_i| \leq 1, \sum y_i = 0 \}$$
taken with respect to $dy_1 \wedge \dots \wedge dy_{d-1}$; certainly each of these volumes is $\leq 1$.

\bibliography{LLLpaper}

\bibliographystyle{amsplain}

  \appendix

   \end{document}